\documentclass{article} 
\newcommand{\ffoot}[1]{}                                 

\newcommand{\bbibitem}{\bibitem}

\newcommand{\llabel}[1]{{\label{#1}}}

\usepackage{graphicx} 

\renewcommand{\r}[1]{(\ref{#1})}
\font\tenmsb=msbm10
\font\sevenmsb=msbm7
\font\fivemsb=msbm5
\newfam\msbfam
\textfont\msbfam=\tenmsb
\scriptfont\msbfam=\sevenmsb
\scriptscriptfont\msbfam=\fivemsb
\def\Bbb#1{{\fam\msbfam\relax#1}}
\newcommand{\bi}{\begin{itemize}}
\newcommand{\ei}{\end{itemize}}
\newcommand{\bd}{\begin{description}}
\newcommand{\ed}{\end{description}}
\newcommand{\be}{\begin{enumerate}}
\newcommand{\ee}{\end{enumerate}}
\renewcommand{\i}{\item}

\newcommand{\bqn}{\begin{eqnarray}}
\newcommand{\eqn}{\end{eqnarray}}
\newcommand{\eqnn}{\nonumber\end{eqnarray}}
\newcommand{\eqnl}[1]{\llabel{#1}\end{eqnarray}}

\newcommand{\nn}{\nonumber}
\newcommand{\noi}{\noindent}
\newcommand{\ba}[1]{\begin{array}{#1}}
\newcommand{\ea}{\end{array}}

\newcommand{\R}{\Bbb{R}}

\newcommand{\fine}{\end{document}}

\def \trait (#1) (#2) (#3){\vrule width #1pt height #2pt depth #3pt}
\def \qed{\hfill
        \trait (0.1) (6) (0)
        \trait (6) (0.1) (0)
        \kern-6pt   
        \trait (6) (6) (-5.9)
        \trait (0.1) (6) (0)
\medskip}
\def \qedmio{\hfill
             \trait (8) (8) (-0.1)
             \medskip}
\def \quadp{{\Huge $\qedmio$}}
\def \quadv{{\Huge $\qed$}}

\newtheorem{ml}{\bf Lemma}
\newtheorem{Theorem}{\bf Theorem}
\newtheorem{mo}{\bf \underline{{\sl Observation}}}
\newtheorem{mrem}{\bf \underline{{\sl Remark}}}
\newtheorem{mcc}{\bf Corollary}
\newtheorem{Definition}{\bf Definition}
\newtheorem{mpr}{\bf Proposition}
\newtheorem{mproperty}{\bf Property}

\newcommand{\bt}{\begin{Theorem}}
\newcommand{\et}{\end{Theorem}}
\newcommand{\bl}{\begin{ml}}
\newcommand{\el}{\end{ml}}
\newcommand{\bo}{\noindent\begin{mo}\rm}
\newcommand{\eo}{\end{mo}}
\newcommand{\bp}{\begin{mpr}}
\newcommand{\ep}{\end{mpr}}
\newcommand{\bc}{\begin{mcc}}
\newcommand{\ec}{\end{mcc}}
\newcommand{\bdeff}{\begin{Definition}}
\newcommand{\edeff}{\end{Definition}}
\newcommand{\bproperty}{\begin{mproperty}}
\newcommand{\eproperty}{\end{mproperty}}
\newcommand{\brem}{\begin{mrem}\rm}
\newcommand{\erem}{\end{mrem}}

\newcommand{\ppotR}[3]
{

\begin{figure}\begin{center}
~\includegraphics[width=#3truecm]{#1.eps}\\
\caption{#2}
\llabel{#1}
\end{center}
\end{figure}
\noindent$\!\!$}
       
\newcommand{\lam}{\lambda}
\newcommand{\ga}{\gamma}\newcommand{\g}{\gamma}
\newcommand{\al}{\alpha}
\newcommand{\eps}{\varepsilon}

\newcommand{\ca}{c_{\alpha}}
\newcommand{\sa}{s_{\alpha}}
\newcommand{\caq}{c_{\alpha}^2}
\newcommand{\saq}{s_{\alpha}^2}

\newcommand{\con}{{\cal C}}

\renewcommand{\H}{{\cal H}}

\newcommand{\cc}{ constant control }  

\newcommand{\da}{\Delta_A^{-1}(0)}
\newcommand{\db}{\Delta_B^{-1}(0)}     
\newcommand{\neigh}{neighborhood \ } 

\newcommand{\beq}{\begin{equation}}     
\newcommand{\eeq}{\end{equation}}     

\newcommand{\adj}{~\mbox{ad}} 
\newcommand{\la}{\lambda}      
\newcommand{\vep}{\varepsilon}
\newcommand{\Si}{\Sigma} 
\newcommand{\vp}{\varphi}

\newcommand{\NH}{N\!\!\!H}
\newcommand{\UNH}{\underline{\NH}}
\newcommand{\NA}{N_{A}}

\begin{document} 

\begin{center} \noindent
{\LARGE{\sl{\bf Time Optimal Synthesis for Left--Invariant Control
Systems  on $SO(3)$}}}
\end{center}

\vskip 1cm
\begin{center}
Ugo Boscain, 

{\footnotesize SISSA-ISAS, 
Via Beirut 2-4, 34014 Trieste, Italy}

Yacine Chitour

{\footnotesize  Universit\'e Paris XI,
D\'epartement de Math\'ematiques, 
F-91405 Orsay, France}
\end{center}

\vspace{.5cm} \noindent \rm 

\begin{quotation}

\noindent  {\bf Abstract}        
Consider the control system $(\Sigma)$ given by $\dot x=x(f+ug)$, 
where $x\in SO(3)$, $|u|\leq 1$ and $f,g\in so(3)$ define two 
perpendicular left-invariant vector fields normalized so that 
$\|f\|=\cos(\al)$ and $\|g\|=\sin(\al)$, $\al\in ]0,\pi/4[$. In 
this paper, we provide an upper bound and a lower bound for $N(\alpha)$, 
the maximum number of switchings for time-optimal trajectories 
of $(\Sigma)$. More precisely, we show that 
$N_S(\al)\leq N(\al)\leq N_S(\al)+4$, where $N_S(\al)$ is a suitable 
integer function 
of $\al$ such that $N_S(\al)\stackrel{\mbox{{\Large $\sim$}}}{\mbox{{\tiny 
$\al\to 
0$}}}\pi/(4\alpha).$
The result is obtained by studying the time optimal
synthesis of a projected control problem on $\R P^2$, 
where the
projection is defined by an appropriate Hopf fibration.
Finally, we study  the projected control problem on the unit sphere $S^2$. 
It exhibits interesting features which will be partly rigorously 
derived and partially described by numerical simulations.
\end{quotation}

\section{Introduction}
Let $(\Sigma)$ be the control system given by:
\beq\llabel{sys0}
\dot x = x(f +ug),
\eeq
where $x\in SO(3)$, $|u|\leq 1$ and $f,g\in so(3)$ give rise to
two non zero perpendicular  left-invariant vector fields on $SO(3)$. 
In this paper, we consider the following problem:
given any pair of points $x_1,x_2$ of $SO(3)$, find a trajectory of
(\ref{sys0}) steering $x_1$ to $x_2$ in minimum time. That issue is
known as the problem of determining the time optimal synthesis (TOS) for  
$(\Sigma)$. The strategy to determine a TOS usually consists in two steps:
\begin{description}
\item[1.] Reduction procedure: it is based on the Pontryagin Maximum, 
Principle (PMP) which is a first-order necessary condition for optimality. 
Roughly speaking, the PMP reduces the candidates for time optimality to 
the so called extremals, which are solutions of a pseudo-Hamiltonian system. 
This reduction procedure may be refined using higher order conditions,
such as Clebsch-Legendre conditions, higher-order maximum principle, envelopes,
conjugate points, index theory (cf. for instance 
\cite{agra-icm,agra-sympl-x,agra-sympl,libro,jurd-book,scha-kre,lm, 
pontlibro,sar,scha,stef,sus1,sus2,
suss-env1,suss-env2});\ffoot{aggiungere poi citare il tedesco e fatima 
leite}
\item[2.] Selection procedure: it consists of selecting the time optimal
trajectories among the extremals that passed the test of Step {\bf 1.} (see 
for instance \cite{bolt,libro,bru,pic-suss}).
\end{description}
Step {\bf 1.} is already non trivial and, in general, the second one is
extremely difficult: if the state space is two-dimensional, the problem of 
determining the TOS for single-input control systems is now 
well-understood 
\cite{ex-syn,automaton,libro,tre,quattro,uno,due,sus1,sus2}. However, for 
higher dimensions, very few examples of complete TOS for a non linear 
control system are available (see for instance \cite{lau-sou}). Intermediate 
issues were thus deeply investigated: determining estimates for the number of 
switchings of optimal trajectories, describing the local structure of 
optimal trajectories, finding families of trajectories sufficient for  
optimality, cf. 
\cite{mario,boschi,scha-kre,mar-pic,pontlibro,scha,suta}, etc.

For the control system $(\Sigma)$, we normalize the two perpendicular
vector fields induced by $f$ and $g$ in such a way that
$\|f\|=\cos(\al)$, $\|g\|=\sin(\al)$, with $\al\in]0,\pi/2[$ (for the
precise meaning of ``perpendicular'' and of the symbol $\|.\|$, we
refer to Section \ref{s-basicfacts}).  Defining $X_+:=f+g$ and
$X_-:=f-g$, we have $\|X_+\|=\|X_-\|=1$ and $\al$ is the angle between
$f$ and $X_+$.

By a standard argument (see Section $2$ below), one can show that
every time optimal trajectory is a finite concatenation of bang arcs
(i.e.  $u\equiv \pm 1$) or singular arcs ($u=0$) and thus, the Fuller
phenomenon (i.e. existence of a trajectory  of a control system joining two
points in (finite) minimum time, with an infinite
number of switchings, cf. \cite{kup,zelibor}) never occurs. (A switching 
time -- or
simply a switching -- along an extremal is a time $t_0$ so that the
control $u$ is not constant in any open neighborhood of $t_0$.)
Moreover, one can easily show that, 
\ffoot{(of course  excepted 
$0$ and $\pi$ where $f$ or $g$ are vanishing and we loose 
controllability)}
the supremum $N(\alpha)$ of the number of switchings 
over all time optimal trajectories of $(\Sigma)$, is finite.

By using the index theory developed by Agrachev,  
it is proved in \cite{agra-sympl-x} that 
\ffoot{CITARE LA FORMULA CON L"INDICE, DIRE DEL MISPRINT, DIRE 
CHE IN REALTA" E" UNA CONDIZIONE SUFFICIENTE DI OTTIMALITA" 
LOCALE}
\begin{equation}\label{agr0} 
N(\alpha)\leq\NA:=\left[\frac{\pi}{\alpha}\right], 
\end{equation} 
where
$\left[\cdot\right]$ stands for the integer part.  That result was not
only an indirect indication that $N(\alpha)$ would tend to $\infty$ as
$\alpha$ tends to zero, but it also provided a hint on
the asymptotic of
$N(\alpha)$ as $\alpha$ tends to zero. 

A related line of work regards 
the study of the distributional version of $(\Sigma)$, which is the
driftless control system given by $\dot x=x(u_1f_1+u_2f_2)$,
$|u_1|,|u_2|\leq 1$ and $f_1,f_2\in so(3)$ linearly independent. Indeed,
assuming that $\|f_1\|=\|f_2\|$, Sussmann and Tang (\cite{suta}) showed that time
optimal trajectories have at most four switchings and they provided a
finitely parametrized family of trajectories sufficient for
optimality.  That result was extended to the general case ($f_1$ and
$f_2$ just linearly independent, cf \cite{boschi}): time optimal
trajectories have at most five switchings.  For both works, the
elimination from optimality of extremals with respectively five or six 
bangs
relies on the envelope theory developed in the context of control theory
by Sussmann (cf. \cite{suss-env2}). 

At the light of the previous results, there was strong evidence for 
two radically situations as $\alpha$ tends to zero: for $(\Sigma)$,
$N(\alpha)$ is expected to go to infinity, as for the distributional
control system, there exists a universal bound on the number of 
switchings. The main result of the present paper confirms that difference,
i.e. $N(\al)$ tends to $\infty$ as $\alpha$ tends 
to zero. More precisely, we complete the inequality (\ref{agr0}) as follows:
\bt\llabel{t-th1}
Let $(\Sigma)$ be the control system defined in (\ref{sys0}) with $f,g$
perpendicular so that $\|f\|=\cos(\al)$ and $\|g\|=\sin(\al)$, 
$\al\in ]0,\pi/4[$. Then, if $N(\alpha)$ is the 
maximum number
of switchings along a time-optimal trajectory of $(\Sigma)$, we have
\beq\label{foo} 
N_S(\alpha)\leq N(\al)\leq 
N_S(\alpha)+4,
\mbox{ where } N_S(\alpha):=2\left[\frac{\pi}{8\alpha}\right]-
\left[2\left[\frac{\pi}{8\alpha}\right]-\frac{\pi}{4\alpha}\right].
\eeq
\et
The above theorem improves (\ref{agr0}) in two ways:
{\bf i)} for $\alpha$ small, it (essentially) divides the upper 
bound of $N(\alpha)$ by  four with respect to (\ref{agr0});
{\bf ii)} it provides a lower bound of $N(\alpha)$ differing from
the upper bound by a constant.

The  lower bound is in fact our main contribution and,
to get it, one must prove the existence of time optimal trajectories
of $(\Sigma)$ admitting {at least} a number of switchings equal to
that lower bound. Our strategy consists of projecting the control
problem onto another $(\Sigma)_S$ defined next. First, let $\R P^2$ be
the two-dimensional real projective space (i.e. the two-dimensional
manifold made of the directions of $\R^3$) and fix a point $x_0\in
SO(3)$. Consider the Hopf fibration $\Pi:SO(3)\rightarrow \R P^2$
defined by $Ker d\Pi(x_0)=Span\{x_0 f \}$, which means, roughly
speaking, that $\Pi$ annihilates the drift term $f$ at $x_0$. Then, we
project $(\Sigma)$ by $\Pi$ and obtain a single-input 
$SO(3)$-equivariant control system $(\Sigma)_S$ on $\R P^2$ given by $\dot
y=y(f_S+ug_S)$, with $f_S=d\Pi(f)$ and $g_S=d\Pi(g)$, {that is locally controllable.}
 We then 
consider
the minimum time problem of connecting $\Pi(x_0)$ to any other point
of $\R P^2$. 

In fact, we study a
slightly different time optimal problem by lifting $(\Sigma)_S$ to the unit
sphere $S^2$. By an abuse
of notation, we still denote by $(\Sigma)_S$ the control system
obtained in that way. Hence $\Pi(x_0)$ is identified with the north pole 
and $\R P^2$ is identified with $\underline{\NH}$, the subset of the 
sphere made of the union of $\NH$, the (open) top
hemisphere of $S^2$, together with half of the equator. 

The time optimal 
problem consists now of
connecting, in minimum time, the north pole with any point of 
$\underline{\NH}$. {Thanks to the suitable choice of the Hopf 
projection and since  
$\al$ belongs to the interval $]0,\pi/4[$, all extremals
of the projected problem are bang-bang (i.e. they are a finite concatenation of trajectories 
corresponding to controls $+1$ or $-1$).}
Let $N_S(\alpha)$ be the supremum of the number of 
switchings
for time optimal trajectories of $(\Sigma)_S$ starting at the north pole
and ending in $\UNH$ (such trajectories
of $(\Sigma)_S$ are actually entirely contained in $\UNH$, see 
Lemma \ref{l-eq-op}).

The use of the Hopf fibration $\Pi$ is motivated by two facts: first,
every time optimal trajectory for the time optimal problem on
$(\Sigma)_S$ staying in $\UNH$ is the projection by $\Pi$ of a time
optimal trajectory for the time optimal problem on $(\Sigma)$ and
thus, $N_S(\alpha)\leq N(\alpha)$. Taking full advantage of the theory
developed in \cite{libro}, we will actually compute \underline{exactly}
$N_S(\alpha)$ as given in (\ref{foo}).  Second, using the fact that
the fiber above $\Pi(x_0)$ is the support of a singular arc (for this 
problem singular arcs are integral curves of the drift $xf$), we show that 
every regular bang-bang trajectory with 
at least $N_S(\alpha)+5$ cannot be optimal and thus, the upper bound.

It is then clear, by now, that the most delicate part of the argument
relies on the exact determination of $N_S(\alpha)$. This is done by
studying the time optimal synthesis, (TOS for short) for the time optimal problem on
$(\Sigma)_S$. 
Such a TOS is usually constructed, following the theory developed in 
\cite{ex-syn,automaton,libro,tre,quattro,uno,due,sus1,sus2},
recursively on the number of extremals arcs, and by checking at each 
step whether they are optimal or not.
For the problem on $\R P^2$, we are not able 
to complete all the steps of the above construction, which would imply
as a byproduct the existence of the TOS. In particular, we cannot show the
optimality of all the 
extremals (i.e. the trajectories candidate for time  
optimality), but, from their study, we can demonstrate enough partial results
in order compute $N_S(\al)$ precisely and thus to conclude the proof of 
Theorem \ref{t-th1}.

The complete time optimal synthesis is then studied numerically (actually 
on the whole $S^2$)  and is showed in the top of Fig. 
\ref{f-critica-unica}. In 
particular, due 
to the  compactness of $S^2$,  
one of the main issues is to understand the singularities developed by the 
minimum time wave front as it approaches to the south pole.
We provide numerical simulations that describe the evolution of the 
extremal front.
{As $\al\to0$, these numerical simulations suggest the emergence 
of three cyclically alternating patterns of optimal synthesis, each of 
them depending on an arithmetic property of $\al$.}


\ffoot{PARLARE DI OTTIMALITA' LOCALE E GLOBALE}
\ffoot{ it is one of the 
few examples of time optimal synthesis on
a compact manifold and so the minimum }
\ffoot{Moreover these numerical simulations show that all singularities 
of the 
synthesis are stable (in the $\con^3)$-topology of vector fields $f_S$ 
and $g_S$).}   
\ffoot{
\bi
\i e' un pendolo deformato
\i doppia elica
\i uno dei pochi esempi di sintesi ottima su un compatto
\i struttura delle curve di switching, chiusura del culo
\i tutte le singolarita' sono generiche
\i contractibility
\ei
}

{The balance of the paper is organized as follows: Section \ref{s-2} 
collects basic facts relative to the time-optimal trajectories of 
$(\Sigma)$; in Section \ref{s-hopf}, the Hopf fibration is described and 
the proof of Theorem~\ref{t-th1} is provided, assuming some facts about 
the time-optimal synthesis of $(\Sigma)_S$, whose arguments are
deferred to the next section.
In particular we use the expression of $N_S(\al)$ and the relation between 
the length of interior bang arcs for the problem on $\R P^2$.
The construction of the time-optimal synthesis of $(\Sigma)_S$ is 
investigated in Section \ref{s-4}, where an exact computation of 
$N_S(\al)$
is established. We conclude the section with two remarks, the first one
explaining the relation between the TOS 
on the sphere and the TOS of a controlled linear pendulum, the second one
establishing a link with an optimal control problem on $SU(2)$.
Finally, in Section \ref{s-5}, we provide the results of 
the numerical 
simulations completing the study of the time optimal synthesis (in 
particular of the 
possible behaviors in a \neigh of
the south pole) and we propose some open problems stated as conjectures.} 
\section{Statement of the Problem and Properties of Optimal 
Trajectories}
\llabel{s-2}
\subsection{Basic Facts}\llabel{s-basicfacts}
In this paper, we consider the control $(\Sigma)$  given by (\ref{sys0}),
where $x\in SO(3)$, $|u|\leq 1$ and $f,g\in so(3)$.
An admissible control $u$ is a measurable function $u: [a,b]\to [-1,1]$,
where $a,b$ depend (in general)
on $u$, cf. \cite{jurd-book}. A trajectory $\gamma$ of 
$(\Sigma)$ is 
an absolutely
continuous curve $\gamma: J\to SO(3)$, where $J=[a,b]$ is a compact 
segment 
of $\R$ such that there exists an admissible control $u$
for which $\dot \gamma(t)=\gamma(t)(f + u(t) g)$
holds a.e. in $J$. 
We then say that $(\gamma,u)$, defined as before, is an
admissible pair for $(\Sigma)$. 
\bdeff
A trajectory $\gamma$ of $(\Sigma)$, defined on $[a,b]$, is {\it time 
optimal} if, for every trajectory $\gamma'$ of  $(\Sigma)$ defined on 
$[a',b']$ with $\g(a)=\g'(a')$ and   $\g(b)=\g'(b')$, we have $b-a
\leq b'-a'$.
\edeff
The Lie algebra $(so(3),[.,.])$ is isomorphic to the Lie algebra 
$(\R^3,\times)$, where $\times$ denotes the vector product. This 
isomorphism is realized by the map:
\bqn\label{phi}
&&\phi_L:so(3)\to \R^3\\
&&\hspace{-2cm}\phi_L\left(\left(\ba{ccc}
0&-c&b\\
c&0&-a\\
-b&a&0\
\ea\right)\right):=
\left(\ba{c}
a\\b\\c
\ea\right),
\nn
\eqn
and provides an inner product on $so(3)$ given by
$<z_1,z_2>:=<\phi_L(z_1),\phi_L(z_2)>,$
where $z_1,z_2\in so(3)$. The symbol $<.,.>$,  in the right--hand side of 
the above equation, 
stands for the Euclidean inner product of $\R^3$.
With this definition, it follows that: 
$<z_1,z_2>:=-\frac12 Tr(z_1z_2).$
In other words, this scalar product is the opposite of  the Killing 
form on $so(3)$.
In the following, $\|z\|:=\sqrt{<z,z>}$ and  $Id$ is the $3\times 3$
identity 
matrix. We will sometimes consider the $2\times 2$ matrix corresponding to
the planar rotation of angle $\beta$ and we use 
$R_{\beta}$ to denote it. \ffoot{vedere se serve}

In this paper, 
we will assume that $f$ and $g$ are perpendicular and 
normalized so that 
$\|f\|=\cos(\al)$ and $\|g\|=\sin(\al)$, $\al\in ]0,\pi/2[$. 
Here, 
we adopt the following notation used throughout the paper,
$\ca:=\cos(\al)$, 
$\caq:=\cos^2(\al)$, $\sa:=\sin(\al)$ 
and $\saq:=\sin^2(\al)$.  We define $h:=[f,g]=fg-gf$ and 
$$
X_+:=f+g,\ \ X_-:=f-g.
$$
Note that $\|X_\eps\|=1$, with $\eps=+,-$. For a 
vector field $z\in so(3)$, we use 
$e^{tz}$ to denote the flow of $z$, acting on the right, so that
$t\mapsto pe^{t z}\in SO(3)$ is the integral curve of $z$ starting at $p$ at
time $0$. Since $z$ is linear, we have 
$e^{tz}=\sum_{n=0}^{\infty}\frac{(tz)^n}{n!}$. We use $\adj_z$ to denote 
the operator $w\mapsto[z,w]$, acting on
vector fields. If $z,w$ are vector fields, then $e^{t\adj
z}(w):=e^{tz}we^{-tz}$. The Lie bracket relations between $f,g,h$ are
$$
[f,g]=h,\ \ [g,h]=\saq f,\ \ [h,f]=\caq g.
$$
{}From them, one deduces the following classical relations that 
will be 
useful later:
\bqn
e^{t\adj X_\eps}(f)&=&(\caq+\saq \cos(t))f+\eps\caq 
(1-\cos(t))g-\eps\sin(t)h,
\llabel{ad1}\\
e^{t\adj X_\eps}(g)&=&\eps\saq (1-\cos(t))f+(\saq +\caq 
\cos(t))g+\sin(t)h,
\llabel{ad2}\\
e^{t\adj X_\eps}(h)&=&\eps\saq 
\sin(t)f-\caq\sin(t)g+\cos(t)h,\llabel{ad3}\\
e^{t\adj X_\eps}(X_{-\eps})&=&(\cos(2\alpha)+2\saq\cos(t))f
+\eps(\cos(2\alpha)-2\caq\cos(t))g-2\eps\sin(t)h,\llabel{ad4}\\
e^{tX_\eps}&=&Id+
\sin(t)X_\eps+(1-\cos(t))X_\eps^2,\llabel{ad5}\\
e^{tadf}(g)&=&\cos(t\ca)g+
\frac{\sin(t\ca)}{\ca}h.\llabel{ad6} 
\eqn
\subsection{Existence of Optimal Trajectories}
A control 
system is {\sl complete} if, for every measurable control 
function $u:[a,b]\to [-1,1]$ and every initial state $p$, there exists a 
trajectory 
$\ga$ corresponding to $u$, which is defined on the whole interval $[a,b]$ 
and satisfies 
$\ga(a)=p$. Since $SO(3)$ is compact and the function ${\cal 
F}(x,u):=x(f+ug)$ is regular enough,
the system \r{sys0} is 
complete.
Note that $(f,g)$ satisfies the Strong Bracket Generating Condition (cf. 
\cite{str}) and the  set of velocities $V(x):=\{x(f+ug),~~u\in[-1,1]\}$ is 
compact and convex. Then, (cf. for 
instance
\cite{suta}): \ffoot{mettere una ref piu' generale}
\bp
\llabel{p-existence-syn}
For each pair of points $p$ and $q$ belonging to $SO(3)$,\ffoot{vedere se 
va bene mettere p e q}
there exists a time optimal trajectory joining $p$ to $q$.
\ep
\subsection{Pontryagin Maximum Principle and Switching Functions}  
\llabel{s-PMP}
We next state the Pontryagin Maximum Principle (PMP) (cf. \cite{pontlibro}) 
for our minimum time problem on $SO(3)$. 
Define the following maps called respectively
\underline{Hamiltonian} and \underline{minimized Hamiltonian}:
\bqn
\H:T^*SO(3)\times [-1,1]\to\R,~~~&\H(p,x,u):=&<p, 
x(f+ug)>,\llabel{inva-ham}\\
H:T^*SO(3)\to \R,~~~&H(p,x):=&\min_{v\in [-1,1]}\H(p,x,v).
\eqnl{brrrr}
The PMP asserts that, if
$\gamma:[a,b]\to SO(3)$ is 
a time optimal trajectory corresponding to a control $u:[a,b]\to [-1,1]$, then 
there exists {\sl a nontrivial field of covectors along} $\ga$, that is an 
absolutely continuous function 
$\la:t\in [a,b]\mapsto \la(t)\in T^\ast_{\ga(t)} SO(3)$ (identified 
with $so(3))$ never
vanishing and a constant $\la_0\geq 0$ such that, for a.e. $t\in Dom(\g)$, 
we have:
\begin{description}
\item[i)]  $\dot \la(t)=-\frac{\partial \H}{\partial 
x}(\la(t),\gamma(t),u(t))=-\lam(t)(f+u(t)g)$,
\item[ii)]  $\H(\ga(t),\la(t),u(t))+\la_0=0$,
\item[iii)] $\H(\ga(t),\la(t),u(t))=H(\ga(t),\la(t))$.
\end{description}
\brem
\llabel{r-postPMP}
The PMP is just a necessary condition for optimality. A trajectory $\g$ 
(resp.  a couple $(\g,\lam)$)  
satisfying the conditions given by the PMP is said to be an {\sl 
extremal} (resp.  an {\sl extremal 
pair}).
An extremal  corresponding to $\la_0=0$ is said to be an  
\underline{abnormal
extremal}, otherwise we call it a \underline{normal extremal}. 
For a normal extremal, we can always normalize $\lam_0=1$, and we do this 
all along the paper.
{Notice that in general an extremal  corresponds to more than one 
covector. For this reason, usually, one distinguishes between abnormal 
extremal that are 
\underline{strict} (i.e. they correspond only to  covectors satisfying 
$\lam_0=0$) and abnormal extremal that are \underline{non-strict} (i.e. 
they correspond to covectors with $\lam_0=0$ and to covectors with 
$\lam_0\neq0$)}. 
\erem
A control
$u:[a,b]\to[-1,1]$ is said to be {\sl bang-bang} if $u(t)\in \{-1,1\}$ 
a.e. in  $[a,b]$. Moreover, if $u(t)\in \{-1,1\}$ and $u(t)$ is constant
for almost
every  $t\in[a,b]$, then $u$ is called a bang control. A {\sl switching} 
time of
$u$ is a time $t\in[a,b]$ such that, for each for every $\vep>0$, $u$ is 
not
bang on $(t-\vep,t+\vep)\cap [a,b] $. A control with a finite number of switchings
is called {\sl regular bang-bang}. A trajectory of $\Si$ is a bang trajectory, 
bang-bang trajectory, regular bang-bang trajectory respectively, if it corresponds 
to a bang control, bang-bang control, regular bang-bang control 
respectively.
The \underline{switching functions}, associated to 
an 
extremal pair $(\g,\lam)$, are the three 
``components'' of the covector $\lam(t)$ on the basis $\{f,g,h\}$ 
transported 
to the point $\g(t)$. More precisely:
\bdeff {\bf (switching functions)}
Let $\Phi_i(x,p)$ $(i=1,2,3)$ be the Hamiltonian functions corresponding 
respectively to the vector fields $f,g,h$ (cf. \cite{jurd-book}). 
i.e. 
$\Phi_1(x,p):= <p,x f>,$
$\Phi_2(x,p):= <p,x g>,$
$\Phi_3(x,p):= <p ,x h>$ and $(\g,\la)$ be an extremal pair. The switching 
functions associated to $(\g,\la)$ are the evaluations of $\Phi_i(x,p)$ 
along the extremal i.e.:
\bqn
\vp_1(t)&:=&\Phi_1(\g(t),\la(t))= <\lambda(t),\g(t) f>,\llabel{vp1}\\
\vp_2(t)&:=&\Phi_2(\g(t),\la(t))= <\lambda(t),\g(t) g>,\llabel{vp2}\\
\vp_3(t)&:=&\Phi_3(\g(t),\la(t))= <\lambda(t),\g(t) h>.\llabel{vp3}
\eqn
\edeff
\brem
\llabel{r-postsw}
Notice that the $\vp_i$'s are at least continuous and
since $\lam$  never vanishes,  the three 
switching functions cannot be all zero at the same time $t$.
Moreover, using the 
switching functions, {\bf ii)} of PMP reads:
\bqn
\H(\la(t),\g(t),u(t))=\vp_1(t)+u(t)\vp_2(t)+\lam_0=0\mbox{ a.e.}
\eqnl{hphi}
\erem
The switching functions are important because they determine 
where the controls may
switch. In fact, using 
the PMP, one easily gets:
\bp A necessary condition for a time
$t$ to be a switching is that $\vp_2(t)=0$. 
Therefore, on any 
interval where 
$\vp_2$ has no zeroes (respectively finitely many zeroes), the
corresponding control is bang (respectively bang-bang). In particular, 
$\vp_2>0$ (resp $\vp_2<0$) on $[a,b]$ implies $u=-1$ (resp.  $u=+1$) 
a.e. on  $[a,b]$. On the other hand, 
if $\vp_2$ has a zero at $t$ and $\dot \vp_2(t)$ exists and is different 
from zero, then $t$ is an isolated switching.
\llabel{p-bangs}
\ep
As a corollary, it holds a.e. along an extremal trajectory that:
\bqn
u(t)\vp_2(t)=-|\vp_2(t)|.
\eqnl{uphi}
An extremal  trajectory $\gamma$ of $\Si$ defined on $[c,d]$
is said to be \underline{singular} if the 
switching function $\vp_2$ vanishes on $[c,d]$. To compute the control 
corresponding to a singular trajectory, one should compute the 
derivatives of the  $\vp_i$'s. Using the Lie bracket 
relations between $f,g,h$, one gets the system of differential equations 
(called the
\underline{adjoint system}) satisfied a.e.:
\bqn
\dot\varphi_1&=& - u \varphi_3,\llabel{p1}\\
\dot\varphi_2 &=& \varphi_3,\llabel{p2}\\
\dot\varphi_3 &=&\saq u \varphi_1 - \caq\varphi_2.\llabel{p3}
\eqn
{}From Eqs. \r{vp3} and \r{p2}, one immediately gets that $\vp_2$ is at 
least a $\con^1$ function.
Moreover, if $\gamma$ is singular in $[a,b]$, then $\vp_2=0$ and, from 
Eq. \r{p2}, we get $\vp_3=0$ a.e.. From \r{hphi} 
(cf. PMP {\bf ii)}), we get that $\vp_1\equiv-1$ a.e. on $[a,b]$. From 
\r{p3} we get $u= 0$ a.e. i.e.:
\bp
\llabel{p-sing}
For the minimum time problem for $(\Sigma)$, singular trajectories are 
integral curves of the drift, i.e. they correspond to a control a.e. 
vanishing. 
\ep  
In the sequel, we will use the following convention. The letter $B$ 
refers to a bang trajectory and the letter $S$ refers to a singular
extremal trajectory. A concatenation of bang and singular trajectories will
be labeled by the corresponding letter sequence, written in order from 
left to right.  
Sometimes, we will use a subscript to indicate the time duration of a trajectory 
so that we use $B_t$ to refer to a bang trajectory defined on an interval 
of length $t$ and, similarly, $S_t$ for a singular trajectory defined on an 
interval of length $t$.

If we fix $u\in[-1,1]$, then the integral curves of $x(f+ug)$ are periodic. 
In particular, the integral curves of $xX_\eps$ are periodic with 
period 
$2\pi$ while the integral curves of the drift $xf$ are periodic with period $2\pi/\ca$.
This means:
\bp
If $\g$ is an extremal trajectory of type $B_t$ (resp. $S_t$), then 
$t<2\pi$ (resp $t<2\pi/\ca$).
\ep
There are two quantities that are remain constant along an extremal 
trajectory. The first one comes from  the fact that the minimized 
Hamiltonian $H$ is 
constant along 
the extremal pairs $(\g,\lam)$ (cf. \r{hphi} and \ref{uphi}):
\bqn
I_1:=-\vp_1(t)+|\vp_2(t)|=\lambda_0,
\eqnl{l0}
with $\lambda_0$ equal to zero or one (cf. Remark \ref{r-postPMP}). 
The second conserved quantity is: 
\bqn
I_2:=\caq\vp_2^2+\saq\vp_1^2+\vp_3^2=K^2, ~~\mbox{for some $K\in\R$.}
\eqnl{el1}

\brem
\llabel{casimiro}
{Equations \r{p1}, \r{p2}, \r{p3} are 
Hamiltonian equations on the dual of 
$so(3)$, with respect to the canonical Poisson structure induced by 
the brackets of  $f,g,h\in so(3)$, and corresponding to the 
left-invariant 
Hamiltonian \r{inva-ham}. The  conserved quantity $I_2$ is the Casimir function 
(see for instance \cite{abhram}).}
\erem

There is a geometric interpretation of the above equations. Let $(\ga,\la)$
be a normal extremal lift of the time-optimal control problem. Then, the
adjoint vector $\la$ with coordinates $(\vp_i)_{i=1,2,3}$, lies in the 
intersection
of the region defined by Eq. (\ref{l0}) and the ellipsoid
defined by Eq. (\ref{el1}).
\ffoot{commentare di piu'} 
\subsection{Classification of optimal trajectories}
{In this section, we investigate the structure of time optimal 
trajectories by  analyzing the extremal flow defined in 
(\ref{p1})-(\ref{p3}), subject to  (\ref{l0}) and (\ref{el1}). 
First we study abnormal extremals (we prove that they are regular 
bang-bang and we establish a relation between the interior bang times). 
Then we study normal extremals that are bang bang (again we find  relation between the 
interior bang times). Finally we study optimal trajectories  
containing a singular arc.} The results presented in this section are 
well-known, and some of them  already
contained in \cite{agra-sympl-x,agra-book}, although in many cases without 
proof. To have a self-contained paper, we provide an argument for all of them.
\subsubsection{Abnormal Extremals}
The following proposition describes the switching behavior of abnormal 
extremals.
\bp
\llabel{p-ab}
Let $\g$ be an abnormal extremal. Then, it is regular bang-bang 
and the time duration between two
consecutive switchings is always equal to $\pi$. In other words, $\g$ is of 
 kind $B_{\pi}B_{\pi}...B_{\pi}B_t$ with $t\leq\pi$.
\ep
{\bf Proof of Proposition \ref{p-ab}}
By definition, $\lambda_0=0$. Then Eq. \ref{hphi} becomes
\bqn
\vp_1(t)=-u(t)\vp_2(t),~\mbox{ for a.e. $t\in Dom(\g)$}.
\eqnl{hphiab}
If $\g$ is singular on some 
interval $[c,d]$, then $\vp_2\equiv 0$ and from \r{p2} $\vp_3\equiv 0$
on $[c,d]$. Eq. \r{hphiab} 
gives $\vp_1\equiv 0$, contradicting
the non triviality of $\lambda$ (cf. Remark \ref{r-postsw}).
 Then $\gamma$ cannot contain a singular arc.
Therefore, $u^2=1$ a.e. $t\in Dom(\gamma)$.

{}From Eqs. \r{p3} and \r{hphiab}, we get a.e. 
$\dot\varphi_3(t)=(-\saq u(t)^2-\caq)\vp_2(t)=-\vp_2(t)$. This means that, 
in the  $(\vp_3,\vp_2)$ plane, the vector 
$z(t):=(\vp_3(t),\vp_2(t))$ 
rotates with angular velocity equal to one (cf. Eq. \r{p2}).
This implies $\gamma$ is a regular 
bang-bang trajectory and the time duration between two 
consecutive switchings along $\gamma$ is always equal to $\pi$. \quadp
\subsubsection{Normal Bang-Bang Extremals}
Let $\gamma$ be a bang-bang trajectory starting at $p_0$ and ending at 
$p_0 e^{(t_0 X_+)} e^{(t_1 X_-)} e^{(t_2 X_+)} e^{(t_3 X_-)}$.
The case in which the first bang is of kind $X_-$ is 
similar.
We have $\vp_2(t_0)=\vp_2(t_0+t_1)=\vp_2(t_0+t_1+t_2)=0$ which implies: 
\beq\label{sw1}
<\la(t_0+t_1),p_2e^{-t_1adX_-}(g)>=<\la(t_0+t_1),p_2g>=
<\la(t_0+t_1),p_2e^{t_2adX_+}(g)>=0, 
\eeq 
where $p_2=p_0 e^{(t_0 X_+)} e^{(t_1 X_-)}$.
We need the following definition.
If $z_1,z_2,z_3$ are (possibly time-varying) vector fields of $SO(3)$, the 
application $q\mapsto q(z_1\wedge z_2\wedge z_3)$ is the {\it field of $3$-vectors} 
associated to the $z_i$'s, where $q(z_1\wedge z_2\wedge z_3)$ is an element of 
$\bigwedge^3 T_qSO(3)$, the $3$-fold exterior power of 
$T_qSO(3)$.\ffoot{dire qualcosa di piu' sulla notazione e sulla 
differenza tra $\times$ and $\wedge$}
We now rewrite Eq.
(\ref{sw1}) by using fields of $3$-vectors. We 
obtain: 
\beq\label{sw11} 
g\wedge e^{-t_1adX_-}(g)\wedge e^{t_2adX_+}(g)=0. 
\eeq 
Thanks to (\ref{ad2}), Eq. (\ref{sw11}) 
is equivalent to $r(t_1,t_2)f\wedge g\wedge h=0$\ffoot{da verificare} for 
an appropriate real-valued function $r$. After computations, we get 
$r(t_1,t_2)=\sin(\frac{t_1-t_2}2)$. This implies that $t_1=t_2=t_3$.

Similarly to what we did in the proof of Proposition \ref{p-ab},
consider now a time optimal trajectory of the form $BB_TB$, where
$B_T$ is a nontrivial interior bang arc associated to a normal
extremal and $T\in ]0,2\pi[$. From \r{p3} and \r{hphi} (with $\lam_0=1$), 
we get $\dot\vp_3=-(\vp_2+\saq u)$. Using \r{p2}, this means that the vector 
$z=(\vp_3,\vp_2+u\saq)^T$ satisfies the differential equation 
$$ 
\dot z=\left(\ba{cc}0&-1\\1&0\ea\right)z,\ \ t\in ]0,T[, 
$$ 
with boundary
conditions (the switching conditions imply $\vp_2(0)=\vp_2(T)=0$) 
$z(0)=(\vp_3(0),u\saq)$ and $z(T)=(\vp_3(T),u\saq)$. 
Using $\vp_2(0)=0$ and the fact that  $u>0$ (resp.  $u<0$) 
implies $0>\dot\vp_2(0)=\vp_3(0)$  (resp. $0<\dot\vp_2(0)=\vp_3(0)$), one 
easily gets 
$
\tan(T/2)={\varphi_3(0)}/{(u~\saq)}<0.
$ 
It follows that $T\in(\pi,2\pi)$. In summary, we proved that 
\bp\label{tdura} 
Let $\gamma$ be a bang-bang normal extremal. Then 
the time duration 
$T$ along an
interior bang arc is the same for all interior bang arcs and verifies
$\pi<T< 2\pi$.  
\ep 
\brem
{}From Propositions \ref{p-ab} and \ref{tdura}, we get that, for an extremal  
bang-bang trajectory (normal or abnormal), the time duration
$T$ along an interior bang arc is the same for all interior bang arcs and 
verifies $\pi\leq T< 2\pi$.
\llabel{r-T}
\erem
\subsubsection{Optimal Trajectories Containing a Singular Arc}
The purpose of this paragraph is to describe the structure of time optimal 
trajectories containing singular arcs.
\bp\llabel{sin}
Let $\gamma$ be a time optimal trajectory containing a singular arc.
Then $\gamma$ is of the type $B_tS_sB_{t'}$, with 
$s\leq\frac{\pi}{\ca}$
if $t>0$ or $t'>0$ and $s<2\frac{\pi}{\ca}$ otherwise.
\ep
{\bf Proof of Proposition \ref{sin}}
Let $\gamma$ be a time optimal trajectory containing a singular arc $S_t$, 
$t>0$. From Proposition \ref{p-sing}, we know that $t<2\pi/\ca$.

Assume now that $\gamma$ contains a singular arc and a nontrivial interior 
bang arc. Then, we may assume that $\gamma$ contains a piece of the type 
$S_sB_t$ or $B_tS_s$ (say the first), with $B_t$ a complete bang arc.
Then we have $\vp_2(s)=\vp_3(s)=\vp_2(s+t)=0$. This translates 
to: $g\wedge h\wedge e^{tadX_{\eps}}(g)=0.$ 
Using \r{ad2}, it implies that $\cos(t)=1$, i.e. $t=2\pi$. This  
contradicts the time optimality of $\gamma$. 
Finally, from \r{ad6}, we get 
$
e^{\frac{\pi}{\ca}adf}(g)=-g.
$
{}From this, we deduce  for $t\geq 0$,
$
e^{\frac{\pi}{\ca}f}e^{tX_{\eps}}=e^{tX_{-\eps}}e^{\frac{\pi}{\ca}f}.
$
Therefore, for $t,s\geq 0$, we have
$
e^{sf}e^{tX_{\eps}}=e^{(s-\frac{\pi}{\ca})f}
e^{tX_{-\eps}}e^{\frac{\pi}{\ca}f}.
$
Then, if $s>\frac{\pi}{\ca}$ and $t>0$ and taking into account what 
precedes,
$e^{sf}e^{tX_{\eps}}$ cannot be optimal. \quadp

\subsection{Uniform Bound on the Number of Switchings for Time Optimal 
Trajectories}
For $\alpha\in ]0,\pi/2[$, let $N(\al)$ be the supremum of the number of 
switchings of any time optimal trajectory on $SO(3)$. Thanks to the left-invariance
of the control system \r{sys0}, we may assume that the supremum is taken over 
any time optimal trajectory starting at $Id$.
In this paragraph, we prove the following:
\bp\label{finite}
For $\alpha\in ]0,\pi/2[$, $N(\al)$ is finite (and thus achieved).
\ep
{\bf Proof of Proposition~\ref{finite}} 
Let us first prove that:
\bd
\i[Claim] every optimal trajectory of $(\Sigma)$  
is a finite concatenation of bang and singular arcs.   
\ed 
Let $\gamma:[a,b]\to SO(3)$ be a time optimal
trajectory of $(\Sigma)$.  Let $S$ be the set of zeroes of $\vp_2$ such
that, if $t\in S$, then $\vp_2$ does not vanish identically in some
neighborhood of $t$. Clearly, $S$ is the set of times $t$ such that
$\gamma(t)$ is the junction of two bang arcs or the the junction of a
singular arc and a bang arc. The conclusion follows if $S$ is
finite. Reasoning by contradiction, $S$ must have a
limit point $\bar{t}$. Moreover $\bar{t}\in S$, otherwise $\vp_2$
would vanish identically in a neighborhood of $\bar{t}$, contradicting
the fact that $\bar{t}$ is a limit point of $S$. Note also that
$\dot\vp_2$ is continuous in an open (in $[a,b]$) neighborhood $N$ of
$\bar{t}$ (see Remark \ref{r-postsw}).  By definition of $S$ and 
$\bar{t}$, 
there exists a
sequence $(t_n)$ in $N$ converging to $\bar{t}$ such that
$\vp_2(t_n)\neq 0$. Pick $n$ large enough, so that, if
$[t'_n,t''_n]$ is the maximal subinterval containing $t_n$
with $\vp_2\neq 0$ on $(t'_n,t''_n)$, then $[t'_n,t''_n]\subset N$. Clearly,
$\vp_2(t'_n)=\vp_2(t''_n)=0$, $\gamma$ is a bang arc on $[t'_n,t''_n]$
and $t''_n-t'_n$ tends to zero as $n$ goes to infinity since $t_n$
tends to $\bar{t}$. But, by Proposition~\ref{tdura}, $t''_n-t'_n\geq
\pi$ for $n$ large enough.  So we reached a contradiction and $S$ is
finite. The claim is proved. \quadv

To finish the proof of Proposition \ref{finite}, it remains to show that 
the (finite) number of switchings for any time optimal trajectory is 
uniformly bounded over $SO(3)$. The argument goes by contradiction:
there would exist, then, a sequence of regular bang-bang time optimal 
trajectories $B_{s_n}B_{T_n}\cdots B_{T_n}B_{t_n}$, where 
$s_n,t_n<2\pi$, $\pi<T_n<2\pi$ and the number of switchings $m_n$ goes to 
infinity as $n$ goes to infinity. Therefore, there exists a sequence of 
points $(x_n)$
of $SO(3)$ such that, the minimum time $\tau_n$ needed to connect $Id$ to $x_n$
by a trajectory of $(\Sigma)$, goes to infinity as $n$ goes to infinity.

To reach a contradiction, it is enough to show that there exists a
time ${\cal{T}}$ so that, for every
point $x\in SO(3)$, there exists a trajectory $\gamma$ of $(\Sigma)$ connecting 
$Id$ to $x$ with $T(\gamma)\leq {\cal{T}}$. By a compactness argument and 
thanks
to the $SO(3)$-invariance of $(\Sigma)$, that would result from the following 
fact: there exists $\bar{t}>0$ and an open neighborhood $U\subset SO(3)$ of 
$Id$ such that every point $x\in U$ can be reached from $Id$ in time less or 
equal than $\bar{t}$. The latter simply results from the facts that $(\Sigma)$ 
has the accessibility property and $e^{tf}$ is periodic.\quadp

\brem
Since the degree of non-holonomy of the distribution generated by $(f,g)$ 
is equal to two, then, by standard controllability arguments, one can quantitatively
relate the size of $U$ and $\bar{t}$ as follows: $U$ contains a ball of radius 
$\frac{C\bar{t}}{\alpha^2}$ for some positive constant $C$. Therefore, 
$N(\alpha)$ can be bounded above by $\frac{C'}{\alpha^2}$, for some positive 
constant $C'$.
\erem

\section{The Hopf Fibration and Proof of Theorem~\ref{t-th1}}
\llabel{s-hopf}
\subsection{The Hopf Projection}
In this paragraph, we describe explicitly the 
Hopf projection from $SO(3)$ to $\R P^2$. This projection provides $SO(3)$ 
with a structure of fiber bundle with base $\R P^2$ and fiber 
$S^1$. In the sequel, we use the identification of $\R P^2$ with 
$S^2\setminus\sim$, where $\sim$ is the antipodal map, that is 
$\R P^2$ is the set of rows $(y_1,y_2,y_3)$, $\sum y_i^2=1$, where
$(y_1,y_2,y_3)\sim (-y_1,-y_2,-y_3)$. In the sequel, $\R P^2$
is identified with the subset $\underline{\NH}$, made of the 
(open) top hemisphere together with half of the equator.
Fix a point $y_0\in\R 
P^2$. The Hopf projection is defined as:
\bqn
&& \Pi:SO(3)\to \R P^2,\nn\\
&& x\mapsto y= y_0 x ,\llabel{hopf1}
\eqnn
where $y_0 x$ is the standard matrix product. Then, any 
left-invariant vector field $V:x\mapsto xv$ on $SO(3)$,
$v\in so(3)$, 
is transformed by $d\Pi$ into the (left-equivariant) vector
field $V_S=d\Pi(V):y\mapsto yv$.
In the following, we call respectively, the control systems $\dot x=x(f+u 
g)$, $x\in SO(3)$ and $\dot y=y (f+u g) $, $y\in \R P^2$, the control 
systems
\underline{upstairs} and \underline{downstairs}.
As explained next, it is crucial to choose $y_0$ so that the drift term 
at the initial point, $Id ~f$, vanishes downstairs, so we require that $\Pi(Id)=y_0$. 
Indeed, if it is the case, then:
{\bf i)} from the point $y_0$, we have 
local controllability (see next section) and this greatly helps in the 
construction of the optimal synthesis;
{\bf ii)} if a trajectory upstairs starts with a singular arc (that is, 
with $u\equiv 0$,
i.e. it is an  integral curve of the drift $f$), then its projection is a point.
That suitable choice of $y_0$ is made possible by the following normalizations:
\bqn
f=\ca \left(\ba{ccc}
0&1 &0\\
-1 &0&0\\
0&0&0
\ea\right) \mbox{ and }
g=\sa \left(\ba{ccc}
0&0 &0\\
0 &0&1\\
0&-1&0
\ea\right),~~~y_0=(0,0,1).
\nn
\eqn
In that way, the control system downstairs reads:
$\dot y=F_S(y)+uG_S(y),$
where
$F_S(y)=y f,$ 
$G_S(y)=y g.$
In order to respect the convention in control theory where states are 
represented by column vectors and costates by row vectors, we will consider
the transposed control system and, with the change of notations $y^T\to 
y$, 
$f^T=-f \to f$, ~$g^T=-g\to g$, we obtain downstairs: 
\bqn
\left\{ 
\ba{rcl}
\dot y&=&F_S(y)+uG_S(y),~~~|u|\leq1,~~~\mbox{ ~~where}\\
F_S(y)&=&f y=f\times y=
\ca \left(\ba{ccc}
0&-1 &0\\
1 &0&0\\
0&0&0
\ea\right) \left(\ba{c} y_1\\y_2\\y_3\ea\right)=
\left(\ba{c} -y_2\\y_1\\0\ea\right),\\
G_S(y)&=&g y=g\times y=
\sa \left(\ba{ccc}
0&0 &0\\
0 &0&-1\\
0&1&0
\ea\right)\left(\ba{c} y_1\\y_2\\y_3\ea\right)=
\left(\ba{c} 0\\-y_3\\y_2\ea\right).
\ea\right.
\llabel{giu}
\eqn

\subsection{Proof of Theorem~\ref{t-th1}}
\llabel{s-proofT1}
For the rest of the paper, we assume $\al\in]0,\pi/4[$.
In this Section we prove Theorem~\ref{t-th1}, using a Lemma describing the 
structure of time optimal trajectories of $(\Sigma)_S$ connecting the 
north pole to any point of $\underline{\NH}$, that will be studied in the 
next section.
\bl
\llabel{l-foundamental}
Consider the control system $(\Sigma)_S$ and the time optimal 
trajectories connecting the north
pole to any point of $\UNH$. Then  
\bd
\item[i)] they are regular bang-bang  with same time durations for the
interior bang arcs, i.e., they are of the type
$B_sB_{v(s)}\cdots B_{v(s)}B_t$, where $s\in[0,\pi]$,
$t\in[0,v(s)]$ and
\begin{equation}\label{vs}
v(s)=\pi+2 \arctan(\frac{s_s}{c_s+\cot^2(\alpha)}).
\end{equation}
\item[ii)]Set $N_0(\alpha):=2\left[\frac{\pi}{8\alpha}\right]$. Then the 
maximum number of switching of these trajectories is 
\begin{equation}\label{N_S}
N_S(\alpha)=N_0(\alpha)-\left[N_0(\alpha)-\frac{\pi}{4\alpha}\right].
\end{equation}
The above formula means that $N_S(\alpha)$ can take the values 
$N_0(\alpha)$, $N_0(\alpha)+1$ or $N_0(\alpha)+2$;
\item[iii)] They are projections, through the Hopf map $\Pi$, of time 
optimal trajectories of 
$(\Sigma)$  starting at $Id$.
\ed		
\el
{\bf Proof of Lemma \ref{l-foundamental}}\\
For the Proof of {\bf i)}, see Proposition \ref{p-bangcompleto} and 
Proposition \ref{rel0}. For the proof of {\bf iii)}, see Lemma 
\ref{l-eq-op} 
and Lemma \ref{ns-n}. For the proof of {\bf ii)} see Proposition 
\ref{p00}.  \quadp\\\\
We now prove separately the two inequalities of Theorem \ref{t-th1}.\\\\
{\bf Proof of the inequality $N_S(\alpha)\leq N(\alpha)$.} {}From 
{\bf iii)} of Lemma \ref{l-foundamental}, every 
time optimal trajectory of $(\Sigma)_S$  connecting the north pole to any 
point of $\UNH$ is the projection by $\Pi$ of a time optimal trajectory of 
$(\Sigma)$ with the same time duration (in particular of the time optimal 
trajectory connecting the two fibers). Therefore, 
$N_S(\alpha)\leq N(\alpha)$. \quadv\\
\ppotR{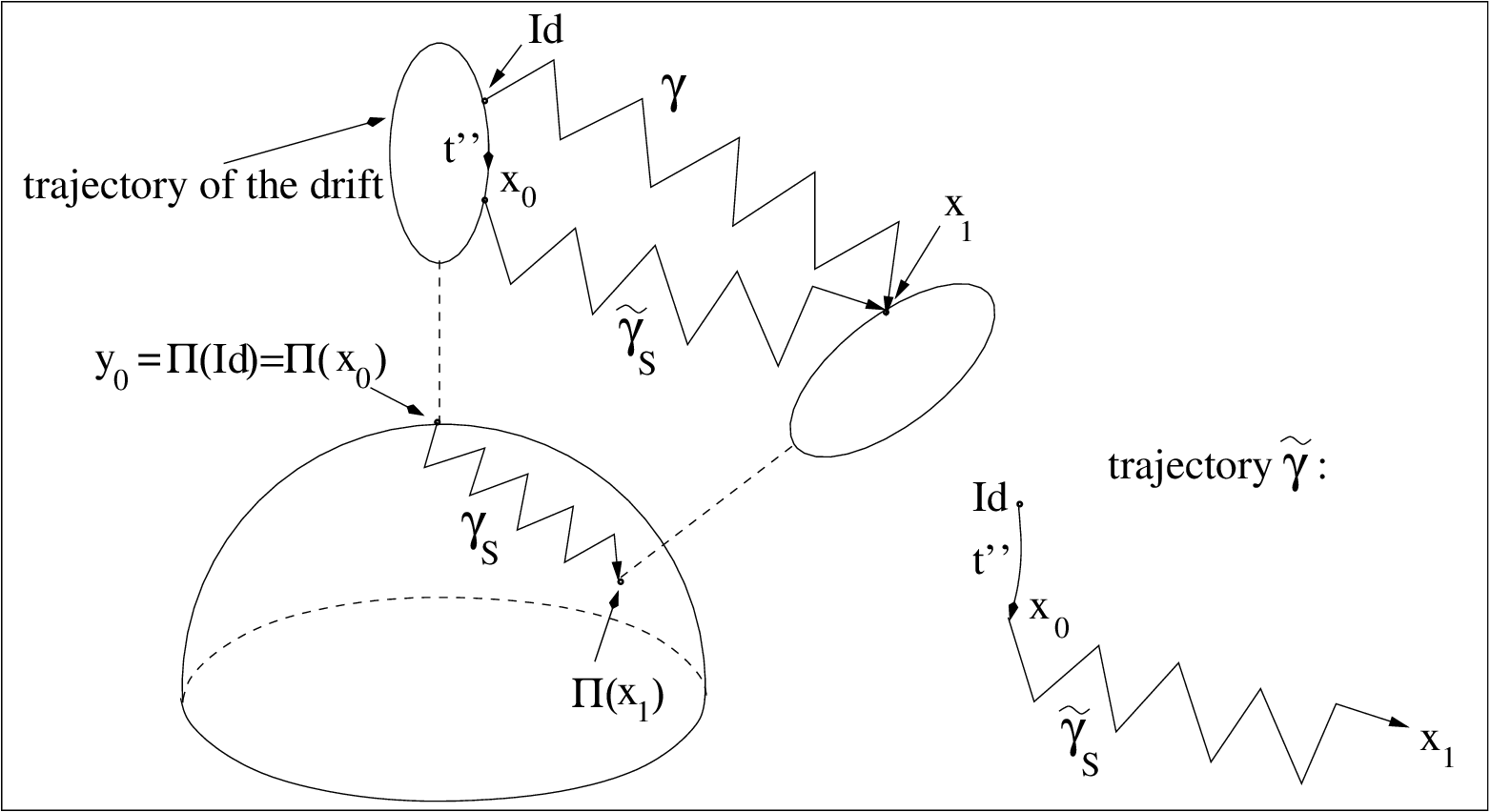}{Proof of Theorem \ref{t-th1}}{12}
{\bf Proof of the inequality $N(\alpha)\leq N_S(\alpha)+4$.} 
We refer to Fig.\ref{f-proof}. Consider a time optimal 
trajectory $\gamma$ of
$(\Sigma)$ containing $N(\alpha)$ switchings. With no loss of
generality, we may assume that $N(\alpha)> 2$.\ffoot{ho messo maggiore e 
non maggiore uguale perche' se c'e' singolre due switching ci possono 
essere} 
By Propositions~\ref{tdura} and \ref{sin}, we deduce that $\gamma$ is regular
bang-bang and is of the type $B_sB_T\cdots B_TB_t$, with $s,t\geq 
0$,\ffoot{forse sarebbe meglio dire che $t\leq T$ e che $s<$??}
$\pi\leq T\leq 2\pi$. Since every subarc of a time optimal trajectory is 
also time optimal, we may assume that $s=t=0$.
Let $Id$ and $x_1$ be the initial and terminal points of $\g$
and consider $\gamma_S$, a time optimal
trajectory for $(\Sigma)_S$ connecting $\Pi(Id)$ and
$\Pi(x_1)$. From {\bf i)} of Lemma \ref{l-foundamental}, $\gamma_S$ is of 
the type $B_{s'}B_{v(s')}\cdots
B_{v(s')}B_{t'}$ with $s'\leq \pi$, $t'< v(s')$ and $m$ interior bangs. 
We thus have $m\leq N_S(\alpha)-1$. We now build, 
from $\gamma_S$,
a suboptimal trajectory connecting $Id$ and $x_1$ as follows: we can
lift $\gamma_S$ to $SO(3)$ to an admissible trajectory
$\tilde{\gamma}_S$ of $(\Sigma)$ connecting $x_0$ and $x_1$, with
$x_0$ in the fiber of $\Pi(Id)$.  It is also clear that $\gamma_S$ and
$\tilde{\gamma}_S$ have same time durations. By construction of the
fiber of $\Pi(Id)$, we get that $x_0=e^{t''f}$ with $t''\leq 2\pi$.
Finally, the curve $\tilde{\gamma}$ obtained as the concatenation of
$e^{t''f}$ and $\tilde{\gamma}_S$ is an admissible trajectory of
$(\Sigma)$ connecting $Id$ and $x_1$. Its time duration is equal to
$$
T(\tilde{\gamma})=t''+T(\tilde{\gamma}_S)=t''+T(\gamma_S)=t''+mv(s')+s'+t',
$$
with $m\leq N_S(\alpha)-1$. Since $\gamma$ is time optimal, we have 
$T(\gamma)\leq T(\tilde{\gamma})$, which implies that
$$
\Big(N(\alpha)-1\Big)T\leq t''+mv(s')+s'+t'.
$$
Using all the estimates on $T,s',t',t''$ (i.e. $T\in[\pi,2\pi)$, cf. 
Remark \ref{r-T},  
$s'\leq\pi,~t'\leq\max_{s'\in[0,\pi]} 
v(s'), ~t''<2\pi$), we deduce that
$$
\Big(N(\alpha)-1\Big)\pi< N_S(\alpha)V(\alpha)+3\pi,
\mbox{ where $V(\alpha):=\max_{s\in [0,\pi]}v(s)$}, 
$$
from which we have
\begin{equation}\label{est11}
N(\alpha)-N_S(\alpha)< N_S(\alpha)\frac{V(\alpha)-\pi}{\pi}+4.
\end{equation}
Set $r(\alpha):=  N_S(\alpha)\frac{V(\alpha)-\pi}{\pi}$. A simple computation
shows that $r(\alpha)=N_S(\alpha)\frac2{\pi}\arcsin(\tan^2(\alpha))$.
Using (\ref{N_S}), it is easy to see that $r(\al)\in]0,1[$ on $]0,\pi/4[$. 
Since $N(\alpha)$ and $N_S(\alpha)$ are integers,
we get $N(\alpha)-N_S(\alpha)\leq 4$. 
\quadv\\\\
\section{The Time Optimal Synthesis Downstairs}\llabel{s-4}

In this section, to compute $N_S(\al)$, we study the time optimal 
synthesis for the problem 
downstairs \r{giu}, starting from the point $y_0$.

\bdeff 
\llabel{d-tos}
{A time optimal synthesis for the problem
downstairs \r{giu}, starting from the point $y_0$ is 
 a family of time optimal trajectories 
$\Gamma=\{\ga_y:[0,b_y]\mapsto \R P^2 $, $y\in\R P^2:~~\g_y(0)=y_0, 
~\g_y(b_y)=y\}$.} 
\edeff

For that purpose, we 
use the theory of optimal syntheses on 2-D manifolds developed  by
Sussmann, Bressan, Piccoli and the first author in 
\cite{ex-syn, morse, automaton,tre,quattro,uno,due,sus1,sus2}  
 and recently rewritten in 
\cite{libro}. 
The  core of the Theory consist of an explicit algorithmic 
construction (by induction on the
number of switchings) of the optimal synthesis.

Note that the previous theory uses a 
a more elaborated concept of synthesis, 
namely   that of \underline{regular synthesis} (see for 
instance \cite{bolt,libro,bru,pic-suss} and cf. Section \ref{s-frames}).
In the following, in order to compute $N_S(\al)$ we just need to follows 
the steps of the algorithmic construction mentioned above, without 
requiring the existence of a regular synthesis. 
In the sequel, by \underline{time optimal synthesis}, we refer to one in 
the sense of Definition \ref{d-tos}, whose existence is simply guaranteed 
by Proposition
\ref{p-existence-syn}.

\bigskip

Consider a two dimensional smooth manifold $M$ and  the problem of 
computing the time optimal synthesis from a fixed point $y_0\in M$ for 
the 
control system: 
\bqn
\dot y=F(y)+uG(y),~~y\in M,~~|u|\leq1,   
\llabel{sys-gen}
\eqn
where                     $F$ and $G$ are  
$\con^\infty$ vector fields.
We introduce three functions:
\bqn
\Delta_A(y)&:=&Det(F(y),G(y))=F_1(y) G_2(y)-F_2(y) 
G_1(y), \llabel{deltaA}\\
\Delta_B(y)&:=&Det(G(y),[F,G](y))=G_1(y)[F,G]_2(y)-G_2(y)[F,G]_1(y),
\llabel{deltaB}\\
f_S(y)&:=&-\Delta_B(y)/\Delta_A(y).\llabel{fs}
\eqn
The sets $\da,\db$ of zeroes of $\Delta_A,\Delta_B$ are respectively
the set of points where $F$ and $G$ are parallel, and the set of
points where $G$ is parallel to $[F,G]$.  These loci are fundamental
in the construction of the optimal synthesis.  In fact, assuming that
they are smooth embedded one dimensional submanifold of $M$ we have
the following:
\bi
\i in each connected region of $M\setminus(\da\cup \db)$, every extremal 
trajectory is bang-bang with at most one switching. Moreover, if the 
trajectory is switching, then the the value of the control switches
from $-1$ to $+1$  if $f_S>0$ and from $+1$ to $-1$ 
if $f_S<0$;\ffoot{dire meglio}

\i the support of singular trajectories (that are trajectories for which 
the switching function identically vanishes, see Definition \ref{d-sw-f} 
below) is always 
contained in the set $\db$;

\i a trajectory not switching on the set of zeroes of $G$ is an abnormal 
extremal (i.e. a trajectory with 
vanishing Hamiltonian) if and only if it switches on the locus $\da$.
\ei
Then the synthesis is built recursively on the number of switchings of 
extremal trajectories, canceling at each step the non optimal 
trajectories (see \cite{libro}, Chapter 1).\ffoot{dire di piu'}

\brem
As we will see later (see Proposition \ref{p-bangcompleto}), the condition 
that $\al<\pi/4$ 
guarantees that 
{\it there are no singular trajectories} 
for the problem downstairs.
\erem

\subsection{Basic Definitions and Facts on Optimal Synthesis on 2-d 
Manifolds}\llabel{sec-bdgiu}

Consider the minimum time problem for the control system \r{sys-gen}.
In this Section, we recall some key facts for the 
construction of time optimal synthesis following \cite{libro}.

The first ingredient is, as usual, the PMP, that, on a two dimensional 
manifold, has exactly the same form as described in Section  
\ref{s-PMP} but with the following change of notation: $x\in SO(3)\to y\in 
M$, $\lam(t)\in T_{\g(t)}SO(3)\to\lam(t)\in T_{\g(t)}M.$
As for the problem upstairs, switchings are described by the switching 
function:
\bdeff {\bf (Switching Function)}
Let $(\g,\lam)$ be an extremal pair. The corresponding
switching function is defined as $\phi(t):=<\lam(t),G(\g(t))>$.
\llabel{d-sw-f}
\edeff
Again, $\phi$ is at least continuously differentiable ($\dot 
\phi(t)=<\lam(t),[F,G](\g(t))>$, cf. discussion in \r{p2},
and it determines the switching rule, according to Proposition 
\ref{p-bangs} with the change of notation $\vp_2\to\phi$. Again, an 
extremal trajectory $\g$, defined on $[a,b]$, is called singular if 
$\phi\equiv0$ in $[a,b]$.
The following three Lemmas illustrate the role of the two functions defined in 
\r{deltaA}, \r{deltaB}. The proofs can be found in 
\cite{ex-syn,libro,due}.
\bl
\llabel{l-sing}
Let $\g$ be an extremal trajectory that is singular in 
$[a,b]\subset Dom(\g)$. Then $\g|_{[a,b]}$ is associated to the so called 
\underline{singular 
control} $\varphi(\g(t))$, where:
\bqn
\llabel{feedbk-varphi}
\varphi(y)=-{\nabla \Delta_B(y)\cdot F(y)\over \nabla\Delta_B(y)\cdot
G(y)},
\eqn
with $\Delta_A$ and $\Delta_B$ defined in Eqs. \r{deltaA} and
\r{deltaB}. Moreover, on $Supp(\g)$, $\varphi(y)$ is always well-defined 
and its absolute value is 
less than or equal to one. 
Finally $Supp(\g|_{[a,b]})\subset \db$.
\el
\bl
\llabel{l-ab}
Let $\g$ be an extremal bang-bang trajectory for the
control problem (\ref{sys-gen}), $t_0\in Dom(\g)$ be a time such that 
$\phi(t_0)=0$ and $G(\gamma(t_0))\neq0$. 
Then, the following conditions are equivalent:
{\bf i)} $\g$ is an abnormal extremal;
{\bf ii)} $\g(t_0)\in\da$;
{\bf iii)} $\g(t)\in\da$, for every time $t\in Dom(\g)$ such that 
$\phi(t)=0$. 
\el     
The following lemma describes what happens when $\Delta_A$ and  
$\Delta_B$ are different from zero. 
\bl
\llabel{l-sw}
Let   $\Omega\subset M$ be an open set such that 
$\Omega\cap (\da\cup \db)=\emptyset$. 
Then all connected components of $Supp(\g)\cap \Omega$, where $\g$ is an 
extremal trajectory of (\ref{sys-gen}), are 
bang-bang with at most one switching. Moreover, if $f_S>0$ throughout 
$\Omega$, then $\g|_{\Omega}$ is associated to a \cc equal to $+1$ or $-1$ 
or has a  switching from $-1$ to $+1$. If $f_S<0$ throughout $\Omega$,
then $\g|_{\Omega}$ is associated to a \cc equal to $+1$ or $-1$ or has a  
switching from $+1$ to $-1$.
\el
\bdeff
\llabel{d-Gamma}
Let $\ga^+:[0,\tau]\to M$ (resp.
$\ga^-:[0,\tau]\to M$) be the trajectory of \r{sys-gen} starting at $y_0$ 
and corresponding to the constant control $u\equiv 1$ (resp. $u\equiv -1$). For 
$t\in]0,\tau]$, let $\Gamma^+(t)$ (resp. $\Gamma^-(t)$) be the support of the curve 
$\ga^+|_{[0,t]}$ (resp. $\ga^-_{[0,t]}$).
\edeff
Under the assumption $F(y_0)=0$ and $\Delta_B(y_0)\neq0$, the next
lemma (for a proof, see for instance \cite{libro,pontlibro}), describes the 
shape of the optimal synthesis
in a \neigh of $y_0$. That local behavior of the optimal synthesis remains
actually the same as long as $\Gamma^+(t)$ and $\Gamma^-(t)$ do not intersect
$\db$ and $\da$ (except of course at $y_0$). 
\ppotR{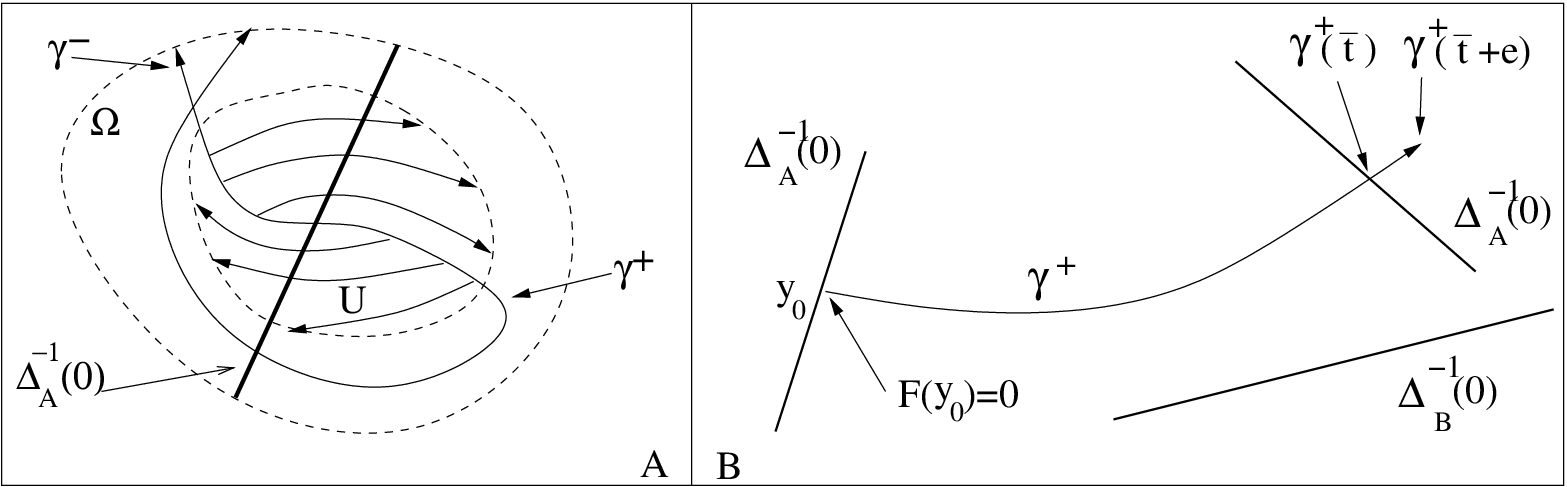}{Lemma \ref{l-origin} and Lemma \ref{l-ultimo}}{14}
\bl 
Consider the control system \r{sys-gen}. 
Assume that $F(y_0)=0$ and $\Delta_B(y_0)\neq0$. 
Let $\Omega$ be an open neighborhood of $y_0$ such that,
$\Omega\cap \db=\emptyset$ and $\Omega \cap\da$ is an embedded 
one-dimensional submanifold of $\Omega$. Let $\ga^+:[0,\tau]\to M$ (resp.
$\ga^-:[0,\tau]\to M$) be the trajectory of \r{sys-gen} starting at $y_0$ 
and corresponding 
to the constant control $u\equiv 1$ (resp. $u\equiv -1$). 
Then, for every $t_+,t_-\in]0,\tau[$
such that: {\bf (a)} $\Gamma^+(t_+),\Gamma^-(t_-)\subset \Omega$,
{\bf (b)} $\Gamma^+(t_+)\cap \da=\Gamma^-(t_-)\cap \da=\{y_0\}$,
{\bf (c)} $\Gamma^+(t_+)\cap \Gamma^-(t_-)=\{y_0\}$, we have the 
following. 
There exists
an open neighborhood $U$ of $\Gamma^+(t_+)\cup \Gamma^-(t_-)$ contained in 
$\Omega$ such that,
for every $y\in U$, there exists a unique extremal trajectory of \r{sys-gen}
of the type $B_sB_t$ contained in $U$, which is time optimal and steers $y_0$ to $y$.
In particular, the system \r{sys-gen} is 
controllable in $U$ and $\ga^+$ (resp. $\ga^-$) is time optimal up to $t_+$
(resp. $t_-$), see Fig.~\ref{fig-mista-1} A.
\llabel{l-origin}
\el
Finally, we need one more lemma, related to Lemma \ref{l-ab}, and 
whose hypothesis are illustrated in Fig.\ref{fig-mista-1} B:
\bl
\llabel{l-ultimo}
Consider the control system \r{sys-gen}. Assume that {\bf i)} $F(y_0)=0$, 
$\Delta_B(y_0)\neq0$, {\bf ii)} there exists $\bar t_+>0$ such that 
$\Gamma^+(\bar t_+)\cap \da =\{y_0,\ga^+(\bar t_+)\}$,
{\bf iii)} there exists $\eps>0$ such that $\Gamma^+(\bar t_++\eps)\cap 
\db=\emptyset$.
Then 
$\g^+$ is extremal exactly up to time $\bar t_+$. Moreover, any extremal trajectory 
$\ga$ defined on $[0,T]$ with $T>\bar t_+$ and coinciding with $\ga^+$ on $[0,\bar t_+]$,
switches at $\bar t_+$ to the constant control $u\equiv -1$ and thus $\ga$ is an 
abnormal extremal (cf. Lemma~\ref{l-ab}).
A similar statement holds for $\ga^-$.
\el

\brem
\llabel{rem-strict}
{Under the hypotheses of Lemma 
\ref{l-ultimo}, one can prove that the abnormal extremal $\g$ 
restricted to an 
interval $[0,\bar T]$,
is a \underline{non-strict abnormal 
extremal} if $\bar T<\bar t^+$, while it becomes a
\underline{strict abnormal extremal} if $\bar T\geq\bar t^+$ (cf. 
Section \ref{s-PMP}). In other 
words, $\g$ becomes a strict abnormal extremal after the first 
switching. These fact are analyzed in details in 
\cite{automaton} and \cite{libro} (see Chapter 
4, and in particular Section 4.3, where strict abnormal extremals are 
called Non Trivial Abnormal Extremals).}
\erem

\subsubsection{Frame Curves and Frame Points}
\llabel{s-frames}

In this paragraph, we briefly recall, for sake of completeness, the main 
results of the theory developed in \cite{uno,due} (see also
\cite{libro}). That material  is only  used here and 
in Section \ref{s-5}, where some numerical
simulations and conjectures are presented.
In \cite{uno,due} (see also \cite{libro}), it was proved that the 
control system \r{sys-gen}, under generic conditions on $F$ and $G$ (with 
the additional assumption $F(y_0)=0$)
admits a time optimal \underline{regular synthesis} in finite time $T$, 
starting 
from $y_0$.
By generic conditions, we mean conditions verified on an open and 
dense subset of the  set of $\con^\infty$ vector fields endowed with 
the $\con^3$ topology (see \cite{libro}, formula 2.6 pp. 39).
 More precisely, let ${\cal R}(T)$ be the 
reachable set in time $T>0$ given by:
\bqn
{\cal R}(T)&:=&\{y\in M :\ \exists ~b_y\in[0,T]\mbox{ and a trajectory
}\nn\\&&\ga_y:[0,b_y]\to M\mbox{ of (\ref{sys-gen}) such that
}\ga_y(0)=y_0,~
\ga_y(b_y)=y\}.
\eqnn
Then a \underline{time optimal regular synthesis} is defined by: {\bf i)} 
a family of time optimal 
trajectories 
$\Gamma=\{\ga_y:[0,b_y]\to M $, $y\in{\cal 
R}(T):~~\g_y(0)=y_0, ~\g_y(b_y)=y\}$ 
such that if $\g_y\in\Gamma$ and $\bar y=\gamma_{y(t)}$ for some $t\in 
[0,b_y]$, then $\gamma_{\bar y}=\gamma_y|_{[0,t]}$; {\bf ii)} a  
stratification of ${\cal R}(T)$ (roughly 
speaking a partition of ${\cal R}(T)$ in manifolds of different 
dimensions, see \cite{libro}, Definition 27, p.56)
such that the optimal trajectories of $\Gamma$ can be 
obtained from a feedback $u(y)$ satisfying:
\bi
\item on  strata of
dimension 2, $u(y)=\pm 1$,
\item  on  strata of dimension 1, called \underline{frame curves}
(FC for short), $u(y)=\pm 1$ or
$u(y)=\varphi(y)$, where $\varphi(y)$ is defined by (\ref{feedbk-varphi}).
\ei
The strata of dimension 0
are called \underline{frame points} (FP).
Every  FP  is an intersection of two FCs. 
\ffoot{c'era: An FP $y$, which is the
intersection
of  two frame curves $F_1$ and $F_2$ is called an
$(F_1,F_2)$ Frame Point.}
 In \cite{due} (see also \cite{libro}), it is
provided a complete classification of all types of FPs and FCs, under
generic conditions. All the
possible FCs are:  
\bi
\i FCs of kind $Y$ (resp. $X$), corresponding to 
subsets
of the trajectories $\g^+$ (resp. $\g^-$) defined as the trajectory 
exiting $y_0$ with \cc $+1$ (resp. \cc $-1$);
\i FCs of kind $C$, called {\sl switching curves}, i.e. curves made of
switching points;
\i FCs of kind $S$, i.e. singular trajectories;
\i FCs of kind $K$,  called overlaps and reached optimally by two
trajectories coming from different directions;
\i FCs  which are  arcs of optimal trajectories starting
at FPs. These trajectories ``transport''  special information. 
\ei
The FCs of kind $Y,C,S,K$ are depicted in Fig.~\ref{fig-FCs}.
There are eighteen topological equivalence classes
of FPs. A detailed description can be found in \cite{morse,libro,due}.
\ppotR{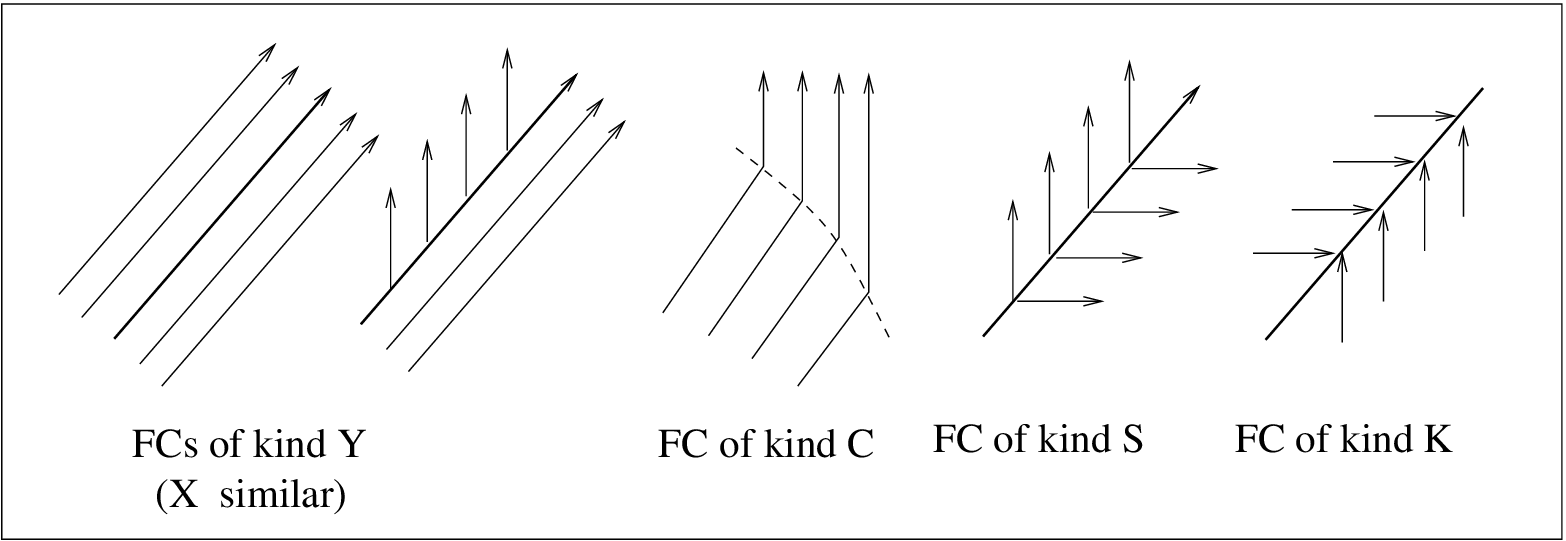}{}{11}
\brem
{The proof of the existence of a regular synthesis is
made  by means of a constructive algorithm (working recursively on the 
number of switchings), that builds  explicitly the 
optimal trajectories (see \cite{libro}, Section 2.5 p.56).
We stress the fact that the existence of a regular synthesis cannot 
be guaranteed before the complete execution of the algorithm.
Since for our systm \r{giu}, we do not reach the end of that
construction, we cannot conclude that such a regular synthesis exists.
However, we conjecture that last fact (see also Section 
\ref{s-5}).
}
\erem
\subsection{The Problem Downstairs}\llabel{sec-PD}

\ffoot{vorrei mettere un lemmino che dica che sotto e sopra abnormal e 
singular sono l'istessa cosa. Per i singular deriva semplicemente dal 
fatto che sia sopra che sotto singular iff u=0. Per gli abnormal bisogna 
oensarci un attimo}

In this section, we apply the theory recalled in Section~\ref{sec-bdgiu} 
to the control system \r{giu} on $S^2$ in order 
to compute $N_S(\al)$, the maximum number of switchings
for time optimal trajectories connecting the north pole to any point of $\UNH$.
First we need some notations.
\bdeff\llabel{long-d}
Set:
\bqn
X^+_S(y)&=&F_S(y)+G_S(y)=X^+y=X^+\times y=
 \left(\ba{ccc}
0&-\ca &0\\
\ca &0&-\sa\\
0&\sa&0
\ea\right) \left(\ba{c} y_1\\y_2\\y_3\ea\right)=
\left(\ba{c} -\ca y_2\\\ca y_1-\sa y_3\\\sa y_2\ea\right),\nn\\
X^-_S(y)&=&F_S(y)-G_S(y)=X^-y=X^-\times y=
 \left(\ba{ccc}
0&-\ca &0\\
\ca &0&\sa\\
0&-\sa&0
\ea\right) \left(\ba{c} y_1\\y_2\\y_3\ea\right)=
\left(\ba{c} -\ca y_2\\\ca y_1+\sa y_3\\-\sa y_2\ea\right).\nn
\eqn
Let $\g:[t_1,t_2]\to S^2 $ be
a  trajectory of (\ref{sys-gen}).
If $\g$ corresponds to the \cc +1 (resp. -1) in $[t_1,t_2]$, we say that
$\g|_{[t_1,t_2]}$ is a $X^+$--trajectory (resp. $X^-$--trajectory).
Moreover, we call $\g^\pm$ the trajectories exiting the point $x_0$ with
respectively, constant control $+1$ and $-1$.  Let $t^\pm_{op}$ be
the last times for
which $\g^\pm$ are optimal. We
define $\g^\pm_{op}:=\g^\pm|_{[0,t^\pm_{op}]}$.
If
$\g_1:[a,b]\to S^2 $ and $\g_2:[b,c]\to S^2 $ are trajectories of
\r{giu} such that
$\g_1(b)=\g_2(b)$, then the {\sl concatenation} $\g_2\ast\g_1$ is the
trajectory:
$$
(\g_2\ast \g_1)(t):=\left\{\ba{l}
\g_1(t)\mbox{ for } t\in[a,b],\\
\g_2(t)\mbox{ for } t\in[b,c].\ea\right.
$$
Notice that, in the notation
$\g_2\ast \g_1$, $\g_1$ comes first.
\edeff
The first quantities to be computed are $\da,\db$ and the sign of $f_S$. 
Referring to Fig.\ref{fig-x+x-}, we have for the system \r{giu}:
\bqn
&&\da=\{(y_1,y_2,y_3)^T\in S^2: y_2=0\},\nn\\
&&\db=\{(y_1,y_2,y_3)^T\in S^2: y_3=0\},\nn\\
&&f_S(y)>0,~\forall~y\in \{(y_1,y_2,y_3)^T\in S^2: y_2y_3>0 \},\nn\\
&&f_S(y)<0,~\forall~y\in \{(y_1,y_2,y_3)^T\in S^2:  y_2y_3<0\}.\nn\\
\eqn
The set $\db$ is called 
the {\sl equator} and $\da$ {\sl the meridian}.
Moreover, let $\NH$ be the (open) top hemisphere, i.e. the set of 
points $(y_1,y_2,y_3)^T$
so that $y_3>0$ and (see Fig.\ref{f-c-1}): 
$$
\NH^+:=\{y\in \NH:y_2<0\},\ \ 
M^+=\{y\in \NH: y_1>0, y_2=0\},\ \ 
E^+=\{y\in S^2: y_2<0, y_3=0\}.
$$
Similarly
$$
\NH^-=\{y\in \NH: y_2>0\},\ \ 
M^-=\{y\in \NH: y_1<0, y_2=0\},\ \ 
E^-=\{y\in S^2: y_2>0, y_3=0\}.
$$
We also parametrize points $y$ of the meridian by the
oriented angle between $\overrightarrow{0y_0}$ and $\overrightarrow{0y}$.  
We use $P(\xi)$,
$\xi\in [-\pi,\pi]$, to denote the point of the meridian
defined by the angle $\xi$. Then $P(0)=y_0$ and $P(\alpha)$ (resp.  
$P(-\alpha)$) is the 
center of rotation in the north hemisphere of $X_S^+$ (resp. $X_S^-$).  We 
also have that $\ga^+$ (resp. $\ga^-$), up to
time $\pi$, is a half-circle with diameter $[y_0,P(2\alpha)]$
(resp. $[y_0,P(-2\alpha)]$), see Fig.~\ref{f-c-1}.
\ppotR{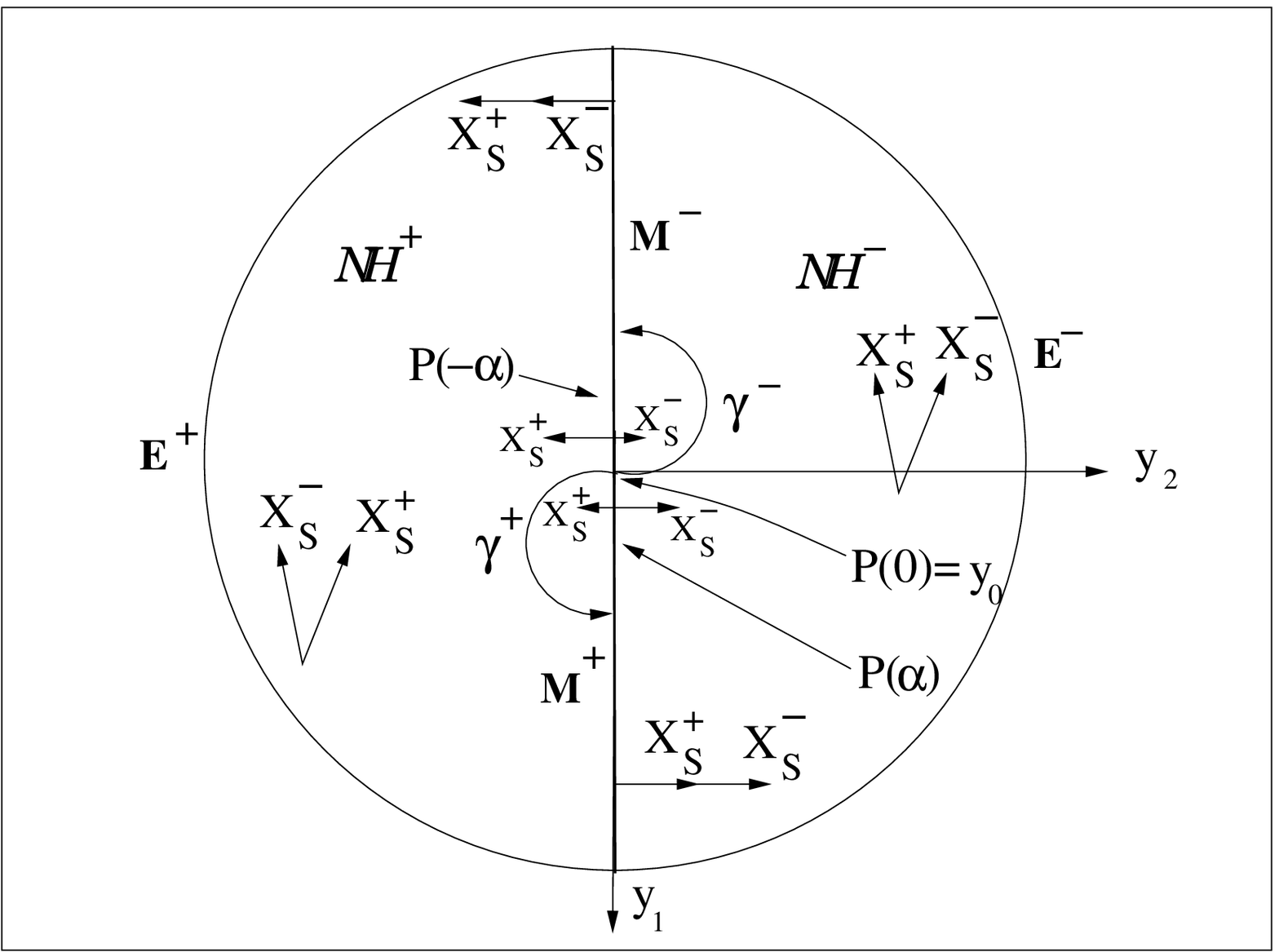}{}{9}
{}From  Lemma
\ref{l-sw}, it follows:
\bp 
Let $\gamma:[0,T]\to S^2$, $\gamma(0)=y_0$ be an optimal trajectory for 
the control system \r{giu}. Then:
\bi
\i $\g$ has at most a $X^+\ast X^-$ switching in $\NH^-$, 
that is, if $Supp(\g|_{[a,b]})\subset  \NH^-$, then 
$\g|_{[a,b]}$ corresponds to one of the three following  controls:
\bd
\i[(-)] $u=+1$ in $[a,b]$,
\i[(-)] $u=-1$ in $[a,b]$,
\i[(-)] there exists $c\in]a,b[,$ such that $u=-1$ in $[a,c[$ and $u=+1$ in $]c,b]$;
\ed
\i $\g$ has at most an $X^-_S\ast X^+_S$ switching in $\NH^+$;
\i $\g$ has at most an $X^-_S\ast X^+_S$ switching in the region $\{x\in 
S^2:y_2>0,y_3<0\}$;
\i $\g$ has at most an $X^+_S\ast X^-_S$ switching in the region $\{x\in 
S^2:y_2<0,y_3<0\}$. 
\ei
\llabel{p-switchings}
\ep
In Fig.~\ref{fig-x+x-}, the integral curves of $F_S,G_S,X^+_S,X^-_S$ and 
the loci $\da,\db$ are depicted. Moreover, the allowed switchings are 
indicated. 
\ppotR{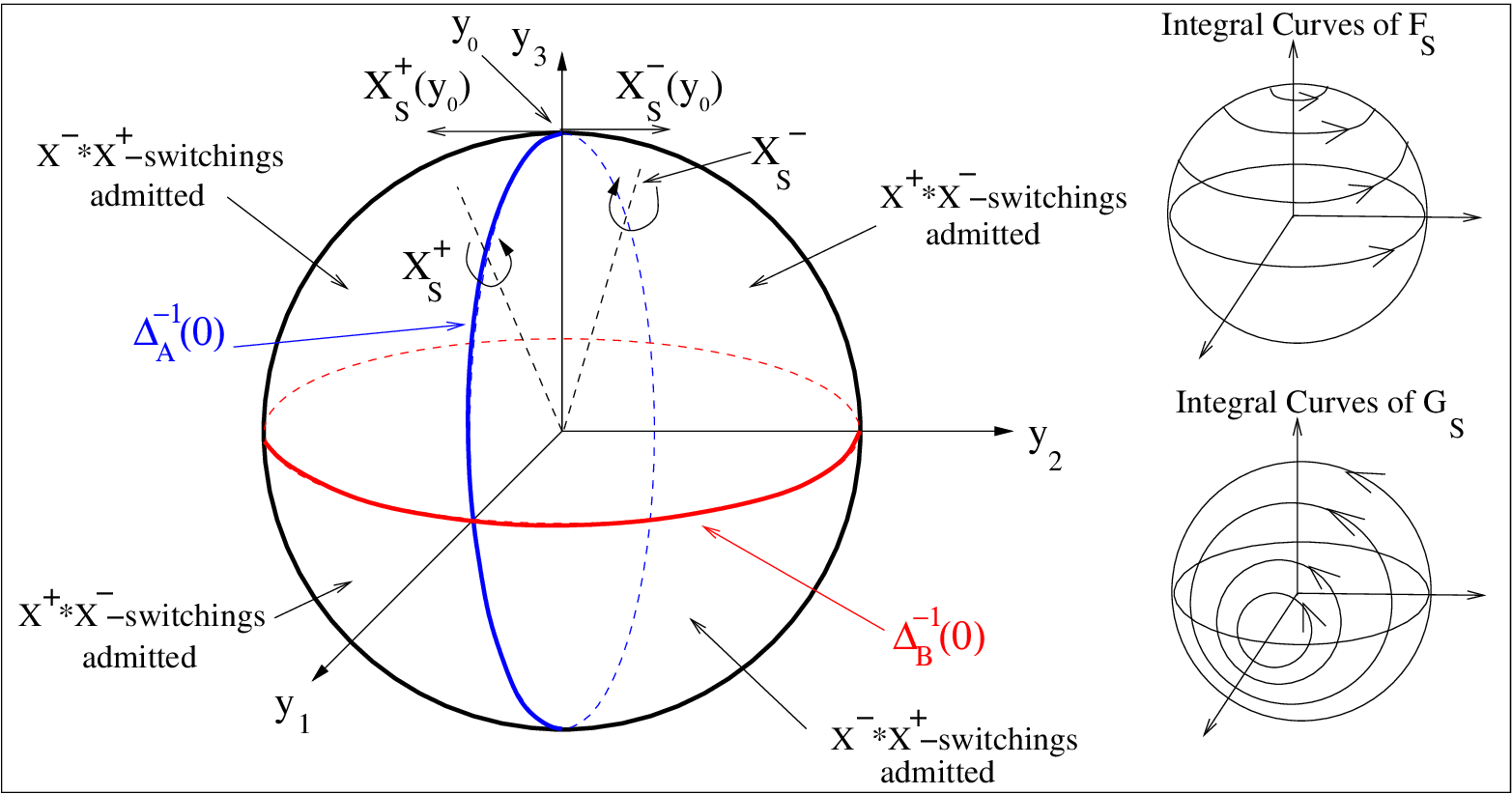}{}{11.5}
\brem\llabel{r-rel-conc}
Notice that, in $\NH^+$ (resp. $\NH^-$), 
$X_S^+$ points on the right (resp.  on the left) of 
$X_S^-$, while, on the meridian, $X_S^+$ and $X_S^-$ are parallel (see 
Fig.\ref{f-c-1}). More precisely, $X_S^+$ and $X_S^-$ point in the 
same direction on $\{P(\xi),\xi\in]\al,\pi-\al[~\bigcup ~]-\pi+\al,-\al[\}$
and in opposite directions on
$\{P(\xi),\xi\in]-\al,\al[~\bigcup~ 
]\pi-\al,\pi]~\bigcup~[-\pi,-\pi+\al[\}$.
\erem

\subsubsection{Two properties of extremal trajectories}
The following two propositions are essential in the construction of the 
optimal synthesis.
\bp
Every time optimal trajectory of (\ref{giu}), starting at the north
pole, is regular bang-bang.  \llabel{p-bangcompleto} 
\ep
{\bf Proof of Proposition~\ref{p-bangcompleto}} Since
$\alpha<\frac{\pi}4$, by taking into account Lemmas~\ref{l-origin} and
\ref{l-ultimo}, the curves $\ga^+$ and $\ga^-$ defined in Definition~\ref{long-d},
do not
intersect the equator and are time optimal until the first time they
meet the meridian, i.e. exactly up to time $\pi$. Moreover, since
singular arcs are contained in the equator and, thanks to
Lemma~\ref{l-ultimo}, any time optimal trajectory $\ga$ of (\ref{giu}),
with at least one switching, is of the form $B_sB_t...$, with $s\in
]0,\pi]$ and $t>0$.  Finally, since $\ga$ is the projection of a time
optimal trajectory $\tilde{\ga}$ of (\ref{sys0}), then the latter is also
of the type $B_sB_t...$. Therefore, by Proposition~\ref{sin}, $\tilde{\ga}$
cannot contain any singular arc, and so $\ga$. \quadp

\bp\label{rel0}
Let $\g:[0,T]\to S^2$ be a time optimal trajectory for the control system 
\r{giu} of the type $B_sB_{t_1}B_{t_2}...$. Then, all time durations of 
interior bang arcs are equal to $v(s)$, where:
\bqn
v(s):=\pi +2\arctan \left(\frac{s_s}{c_s +\cot^2 (\alpha )}\right). 
\llabel{v(s)}
\eqn
\ep
{\bf Proof of Proposition~\ref{rel0}} Consider $\tilde\gamma:[0,T]\to 
SO(3)$, an optimal 
trajectory that projects on $\g$ through the Hopf fibration $\Pi$. Thanks to 
Proposition~\ref{tdura} (see also Remark \ref{r-T}), 
we have $\tilde\g=B_sB_{t_1}B_{t_1}\cdots$, where $t_1\in [\pi,2\pi[$. 
Moreover, since that curve 
projects on 
a time optimal trajectory for \r{giu}, we will establish a relation between $s$ and 
$t_1$. We start from the relations $\vp_2(s)=\vp_2(s+t_1)=0$, which can be written
\beq\label{swisp}
<\la(s),G_S(\g(s))>=<\la(s+t_1),G_S(\g(s+t_1))>=0.
\eeq
Recall that $\la(s)=\la(0)e^{-sX_{\eps}}$, 
$\la(s+t_1)=\la(0)e^{-sX_{\eps}}e^{-t_1X_{-\eps}}$
and $\g(s)=e^{sX_{\eps}}\g(0)$, $\g(s+t_1)=e^{t_1X_{-\eps}}e^{sX_{\eps}}\g(0)$. Since
$\g$ is nontrivial, then $\la(0)$ is a nonzero line vector of $\R^3$. Moreover, 
$\g(0)=y_0=(0,0,1)^T$.
Eq. (\ref{swisp}) can be written as:
$$
\la(0)e^{-sX_{\eps}}(g\times e^{sX_{\eps}}\g(0))=0,\ \ 
\la(0)e^{-sX_{\eps}}e^{-t_1X_{-\eps}}(g\times 
e^{t_1X_{-\eps}}e^{sX_{\eps}}\g(0))=0.
$$
The previous equations can be transformed to:
$$
\det(e^{sX_{\eps}}\la(0)^T,g,e^{sX_{\eps}}\g(0))=0,\ \
\det(e^{t_1X_{-\eps}}e^{sX_{\eps}}\la(0)^T,g,e^{t_1X_{-\eps}}e^{sX_{\eps}}\g(0))=0,
$$
and then to:
$$
\det(e^{sX_{\eps}}\la(0)^T,g,e^{sX_{\eps}}\g(0))=0,\ \
\det(e^{sX_{\eps}}\la(0)^T,e^{-t_1X_{-\eps}}g,e^{sX_{\eps}}\g(0))=0.
$$
Since $e^{sX_{\eps}}\la(0)^T$ is not zero, we deduce that:
\beq\label{sp0}
\det(g,e^{sX_{\eps}}\g(0),e^{-t_1X_{-\eps}}g)=0.
\eeq
We end up with the relation:\ffoot{verificare}
\beq\label{req}
-\sa^2\cos(s-t_1/2)=\ca^2\cos(t_1/2).
\eeq
Taking into account that $\pi\leq t_1<2\pi$, 
we can simplify the previous equation to get
(\ref{v(s)}).\ffoot{METTERE UN POCO DI STUDIO DELLA $v(s)$} \quadp

\subsection{Construction of the time optimal synthesis} 
\llabel{s-ctos}
In this section, we present, step by step, the construction of the TOS
for (\ref{giu}). Since we will not complete that construction, we only 
provide here the steps for which the outcome is justified by a rigorous argument. 
For the other steps of the construction, we refer to the last section where we propose 
conjectures on their outcomes, which are supported by numerical simulation.

\bd
\i[Step $1$] By  Lemmas~\ref{l-origin} and \ref{l-ultimo},  for every
$\eps>0$, there exists an open neighborhood $U$ of $\Gamma^+(\pi-\eps)\cup 
\Gamma^-(\pi-\eps)$ (recall Definition \ref{d-Gamma}) where the time 
optimal synthesis is described
in Fig.~\ref{f-mista-3} A. Moreover, $t^+_{op}=t^-_{op}=\pi$ (recall 
Definition 
\ref{long-d});

\i[Step $2$] {taking into account the analysis of Sections~\ref{sec-bdgiu} and
\ref{sec-PD}, the time optimal trajectories for the problem downstairs are described by the 
following:
\ppotR{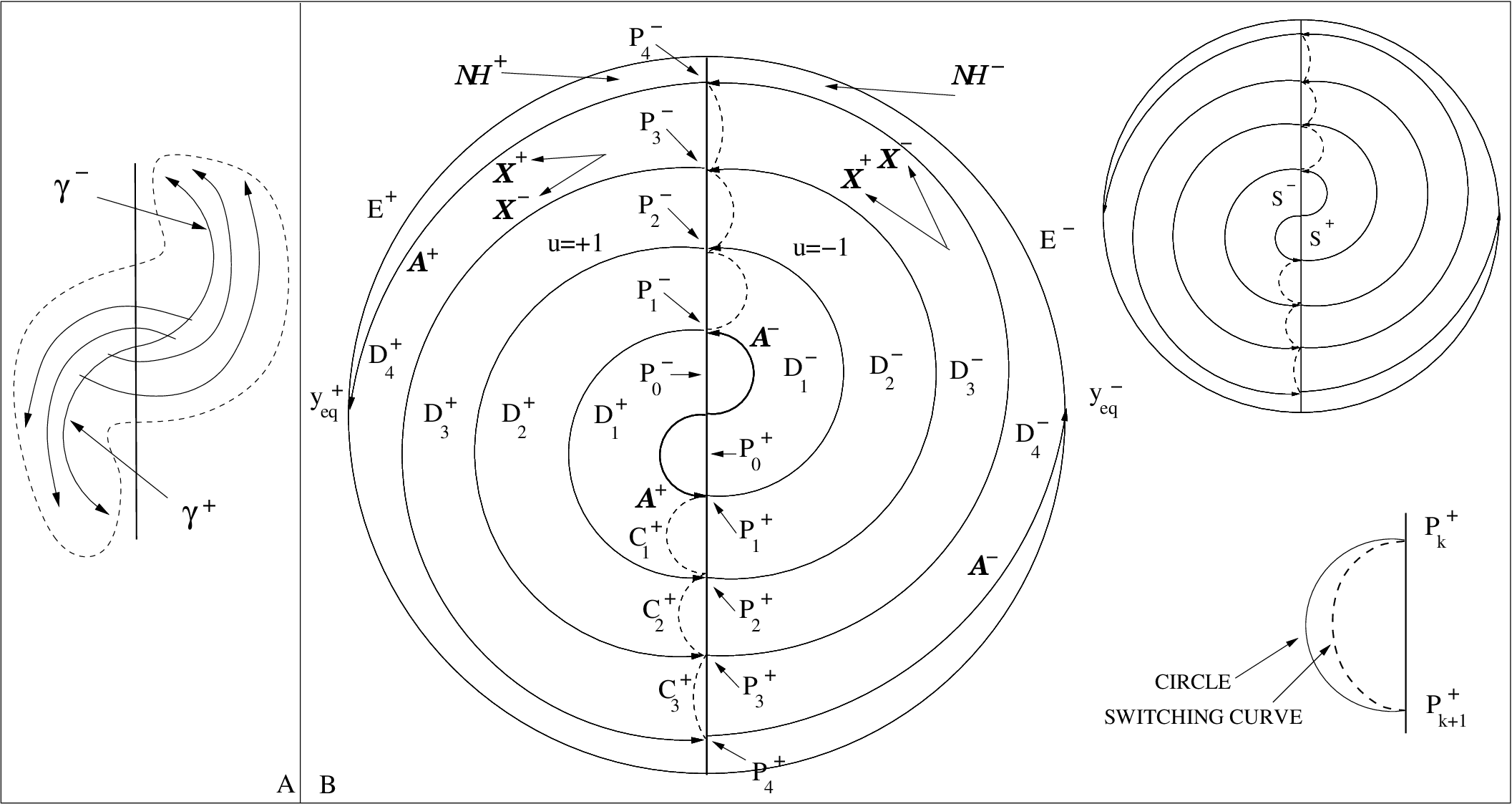}{}{16}
\bp
\llabel{p-sets}
Every time optimal trajectory for the system \r{giu}, starting from the north pole, is 
contained in the following two sets of extremals, 
which are parametrized by the length of the first
bang arc, the one of the last bang arc and the number of arcs:
\bqn
\Xi^+(s,t)&=&\overbrace{e^{X^{\eps}_S t}e^{X^{-\eps}_S v(s)}\cdots 
e^{X^-_S v(s)}e^{X^+_S s}}^{m ~~\mbox{terms}}y_0,\llabel{setEX-1}\\
\Xi^-(s,t)&=&\underbrace{e^{X^{\eps'}_S t}e^{X^{-\eps'}_S v(s)}\cdots 
e^{X^+_S v(s)}e^{X^-_S s}}_{m'~~\mbox{terms}}y_0,\llabel{setEX-2}
\eqn
where $s\in [0,\pi]$, $t\in [0,v(s)]$, the number of bang arcs  ($m$ and 
$m'$ respectively) is an integer and
\bd
\i[(-)] $\eps=+1$ (resp.   $\eps=-1$), if $m$ is odd 
(resp. even),
\i[(-)] $\eps'=+1$ (resp.   $\eps'=-1$), if  $m'$ is even
(resp. odd).
\ed
\ep
}
\i[Step $3$] Let ${\cal{A}}^+$ and ${\cal{A}}^-$ be the two extremal 
trajectories 
starting resp. with controls $u\equiv 1$ and $u\equiv -1$, and switching 
after time $\pi$, i.e. corresponding resp. to
$\Xi^+(\pi,\cdot)$ and $\Xi^-(\pi,\cdot)$. These two curves are abnormal extremals
and their respective first bang arcs coincide with $\ga^+_{op}$ and $\ga^-_{op}$.
{As explained in Remark \ref{rem-strict},  these two curves become 
strict abnormal extremals after time $\pi$.} 

To describe them, consider, for $\eps=\pm$ and $0\leq k\leq \tilde k$ 
($\tilde k$ defined below), the 
half-circles $L^{\eps}_k\subset Clos(\NH^{\eps})$, whose centers lie on 
$\overrightarrow{0P(\eps\al)}$  and passing through the points 
$P^{-\eps}_k$ and $P^{\eps}_{k+1}$, where
$$
P^+_n:=P(2n\alpha),\ \ P^-_n:=P(-2n\alpha),
$$
for the integers $n$ so that $2n\alpha\leq \frac{\pi}2+2\alpha$. Note
that $\frac{\pi}2+2\alpha<\pi$ for $\alpha<\frac{\pi}4$ and, in fact,
the last $P^{\eps}_n$ belongs to the bottom-half hemisphere, i.e. $n\leq 
\tilde k$ where
$\tilde k:=2+\left[\frac{\pi}{4\alpha}\right]$.  It is easy to see that
${\cal{A}}^+$ intersects the top half-meridian according to the
following ordered sequence of points: $y_0$, $P^+_1$, $P^-_2$,
$P^+_3,...$. Similarly, ${\cal{A}}^-$ intersects the
top half-meridian at $y_0$, $P^-_1$, $P^+_2$, $P^-_3,...$.  Moreover,
let $y^+_{eq}$ and $y^-_{eq}$ be the antipodal points of the equator
which are the respective first intersections of ${\cal{A}}^+$ and
${\cal{A}}^-$ with the equator. Note that they are reached at the same
time $T_{eq}$.  Finally, consider the open subset of the
top-hemisphere bounded below by the equator and obtained by removing the
supports of ${\cal{A}}^+$ and ${\cal{A}}^-$ up to time $T_{eq}$,
i.e. all the $L^{\eps}_k$. That set is the disconnected union of the
two ``snake-shaped'' simply connected regions $S^+$ and $S^-$ (defined so
that each $S^{\eps}$ contains the center of rotation of $X_S^{\eps}$). 
Clearly $S^{+}$ and $S^{-}$
are made of open segments of the meridian and open simply connected
regions $D^{\eps}_k\subset \NH^{\eps}$ defined as follows.
For $k=0$, $D^{\eps}_0$ is delimited by $Supp(\ga_{op}^{\eps})$ and the
segment $[P_0,P_1^{\eps}]$ and, for $k\geq 1$, $D^{\eps}_k$ is
delimited by $L^{\eps}_{k-1}$ on the top, $L^{\eps}_k$ on the bottoms
and by the segments $[P^{-\eps}_{k-1},P^{-\eps}_k]$ and
$[P^{\eps}_k,P^{\eps}_{k+1}]$ on the sides, see Fig.~\ref{f-mista-3} B.

In the sequel, if $A,B$ are two subsets of points of $S^{\eps}$, we say
that $A$ is {\it above} $B$ (or equivalently $B$ is {\it below} $A$)
if $A\subset D^{\eps'}_k$ and $B\subset D^{\eps''}_{k'}$ with $k<k'$, for 
some $\eps',\eps''$.

\i[Step $4$] The switching curves (SC for short), associated to the set of 
extremals given in (\ref{setEX-1}) and (\ref{setEX-2}), are defined as 
follows: they can be divided in two families, $(C^+_k)$ and $(C^-_k)$.
If $\eps=\pm$, $1\leq k\leq N_0-1$ and $s\in [0,\pi]$,  then 
$$
C_1^{\eps}(s)= e^{X^{\eps}_Sv(s)}   e^{X^{-\eps}_Ss}  y_0,\ \ \
C^{\eps}_{k+1}(s)=e^{X^{\eps}_Sv(s)}C^{-\eps}_k(s).
$$
The boundary points of $C^{\eps}_k$ are $C^{\eps}_k(0)=P^{\eps}_k$ and
$C^{\eps}_k(\pi)=P^{\eps}_{k+1}$. By using
Proposition~\ref{p-switchings} and since $v(s)\geq \pi$, $s\in
[0,\pi]$, the support of $C^{\eps}_k$ is contained in the subset 
of $Clos(\NH^{\eps})$, delimited by the half-circle 
centered on $\overrightarrow{0P(\al(2k+1))}$ and passing 
through the points $P^{\eps}_k,P^{\eps}_{k+1}$,
and the segment of the
meridian $[P^{\eps}_k,P^{\eps}_{k+1}]$, see Fig~\ref{f-mista-3} B. In
particular, a (SC) with boundary points in the top-hemisphere is
entirely contained in the top-hemisphere and the intersection of its
support with the top-meridian reduces to its boundary points (see 
Lemma \ref{l-ab}).

We next describe the shape of the first (SC) intersecting the equator.
By symmetry, we may assume $\eps=+$. We claim that its intersection with 
the equator
reduces to the point $P(\frac{\pi}2)=(1,0,0)^T$. Indeed, by the
switching rules established in Proposition~\ref{p-switchings}, the
(SC) intersecting the equator is contained in 
$\{y\in S^2: y_2\leq
0,y_3\geq 0\}\cup \{y\in S^2:y_2\geq 0,y_3\leq 0\}$.  Taking into
account the regularity of the (SC) and the values of its boundary
points, the claim is proved, see Fig.~\ref{f-mista-2} A.
\ppotR{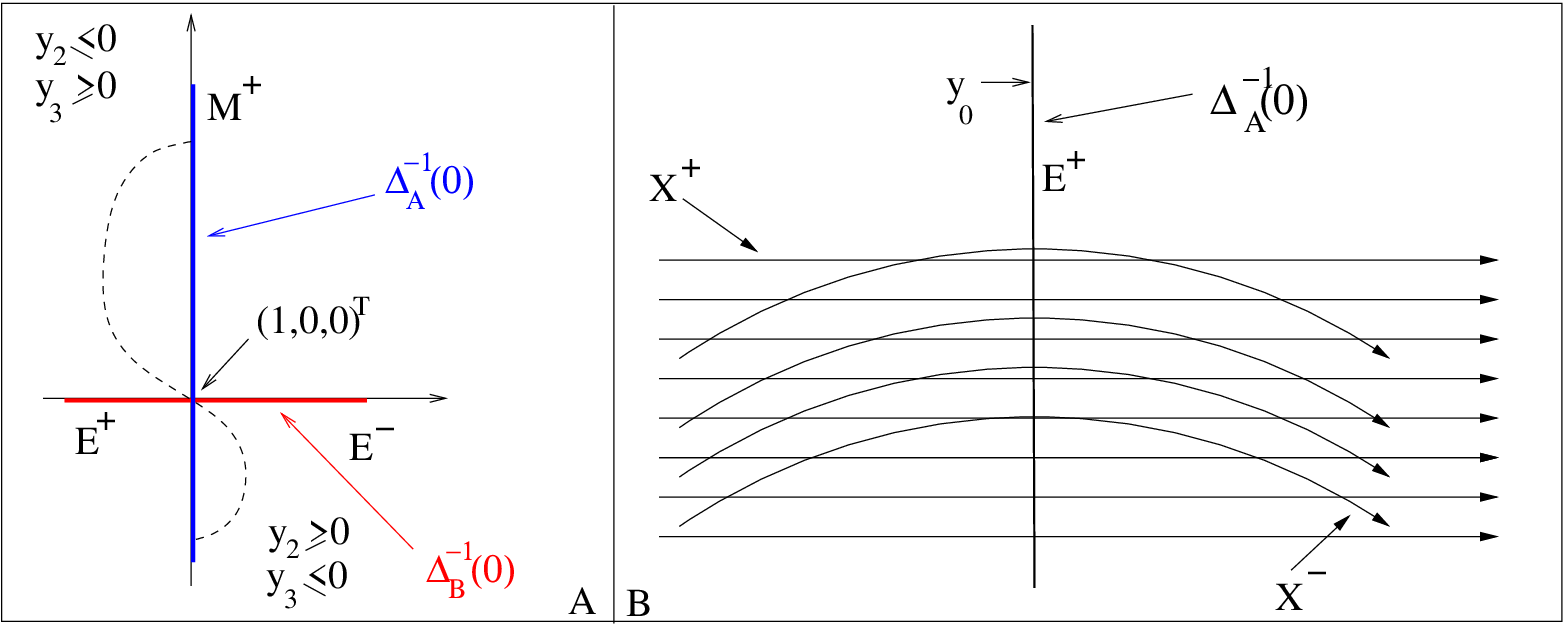}{}{13}
\ed

\subsection{Computation of $N_S(\alpha)$} 
In the previous section, we provided detailed informations about
extremal trajectories and switching curves but we did not show that
every extremal of (\ref{setEX-1}) and (\ref{setEX-2}) is in fact time
optimal.  Anyway, a rigorous derivation of $N_S(\alpha)$ is possible
with the available knowledge of time optimal trajectories combined
with the subsequent lemmas.

\bl\llabel{l-eq-op} 
Every time optimal trajectory $\ga$ starting at $y_0$ 
intersects the equator at most once.
\el
{\bf Proof of Lemma~\ref{l-eq-op}} we argue by contradiction. There would exist
two distinct points of the equator $q_i,q_f$ so that $\ga(t_i)=q_i$,
$\ga(t_f)=q_f$ and $\ga|_{(t_i,t_f)}$ is entirely contained in the (closed)
bottom hemisphere. Let $\ga_{sing}$ be the integral curve of $F_S$ (contained 
in the equator) connecting $q_i$ to $q_f$. Consider now the region of the 
bottom hemisphere bounded by $\ga_{sing}$ 
and $\ga|_{(t_i,t_f)}$. Taking into account, first, the relative positions of 
$X^+_S,X_S^-$, $F_S$ and $G_S$ along the equator and, second, the sign of $f_S$
in the bottom hemisphere, one can check that 
$T(\ga_{sing})\leq T(\ga|_{(t_i,t_f)})$. The argument is similar to that of 
\cite{sus1} (see also \cite{libro}) and is based on the use of Stokes theorem.
Since time optimal trajectories starting at $y_0$ do not contain a singular arc,
it follows that $\ga$ cannot be time optimal. We reached a contradiction.
\quadp

\bl\llabel{ns-n} 
Every time optimal trajectory $\ga$, starting at $y_0$
and remaining in $\UNH$, is the projection of a time optimal trajectory
of $(\Sigma)$ starting at $Id$. \ffoot{questo lemma e' troppo ovvio e poi 
e' gia stato usato per trovare il lower bound. Quindi andrebbe eliminato} 
\ffoot{BISOGNA ANCHE DIRE CHE LE CURVE DI SWITCHING SONO TG AGLI ABNORMAL 
EXTREMALS}

\el 
{\bf Proof of Lemma~\ref{ns-n}} From the definition of the
Hopf fibration, every trajectory $\ga$ of $(\Sigma)_S$, starting at
$y_0$, associated to an admissible control $u$ and staying in
$\UNH$, is the projection of the trajectory $\bar{\ga}$ of $(\Sigma)$
starting at $Id$ with the same control $u$. In particular, $\ga$ and $\bar{\ga}$ 
have same time duration. It is clear that, if $\ga$ is time optimal, then
$\bar{\ga}$ is also time optimal.
\quadp

\bl\llabel{l-snake-op} 
Recall that $S^\eps\subset\NH$.
With the notations
above, pick any point $y$ in the region $S^{\eps}$ and let $\ga_y$ be
a time optimal trajectory connecting the north pole $y_0$ to $y$.  If
$s\in ]0,\pi[$ is the time duration of the first bang arc and $T(y)$
the total time duration of $\ga_y$, then $\ga_{y}|_{(s,T(y)]}$ is
entirely contained in $S^{\eps}$.  
\el 
{\bf Proof of Lemma~\ref{l-snake-op}}
By the switching rules of Proposition~\ref{p-switchings}, along every
time optimal trajectory contained in $\NH^{\eps}$, the control must
switch from $\eps$ to $-\eps$, when arriving at a switching curve
$C^{\eps}_k$. In addition, the time optimal trajectory switches from being an
arc of circle (integral curve of $X^{\eps}_S)$ to another arc of
circle of bigger radius (integral curve of $X^{-\eps}_S)$.  After
rectification of the flow of $X_S^{\eps}$, (i.e. the one entering the 
(SC) $C^{\eps}_k$),  then, by taking into account 
Remark~\ref{r-rel-conc}, 
one gets the situation depicted in Fig.~\ref{f-mista-2} B.

By contradiction, we assume that there exists a time optimal trajectory $\ga$
with time duration $T$ and first bang arc time duration $s<T$ such that
$\ga$ connects $y_0$ to $y\in S^{\eps}$ and $\ga|_{(s,T]}$ exits from $S^{\eps}$.
Let $t'$ be the smallest time (in $[0,T]$) so that $\ga|_{(t',T]}$ is entirely
contained in $S^{\eps}$. 
Then $\ga(t')$ belongs to $Supp({\cal{A}}^+_{[0,T_{eq}]})\cup 
Supp({\cal{A}}^-_{[0,T_{eq}]})$ (see step 3 of Section \ref{s-ctos} for 
the 
definition  of $T_{eq}$).
If $\ga(t')$ is on the (top)-meridian, then it has to switch so that the interior
bang time duration is constant, equal to $\pi$. Therefore $\ga|_{(t',T]}$
will never re-enter $S^+$. We thus deduce that $\ga(t')$ is not on the meridian and,
with no loss of generality, we will assume that $\ga(t')$ belongs to the 
($1$-dim.) 
interior of some $L^+_k$, $k\geq 1$. Now we make the following two claims:\\\\
{\bf Claim 1}
{\it With the notations above, there exist $t''<t'<t'''$ such that
$$
\ga|_{(t',t''')}\subset D^{\eps'}_k\subset S^{\eps}\mbox{ and }
\ga|_{(t'',t')}\subset D^{\eps'}_{k+1}\subset S^{-\eps},\mbox{ for some } 
\eps'\in\{+,-\},
$$
i.e. $\ga$ passes (backward in time) from $S^{\eps}$ to 
$S^{-\eps}$ at time $t'$ by going ``down''.
}\\\\
\ppotR{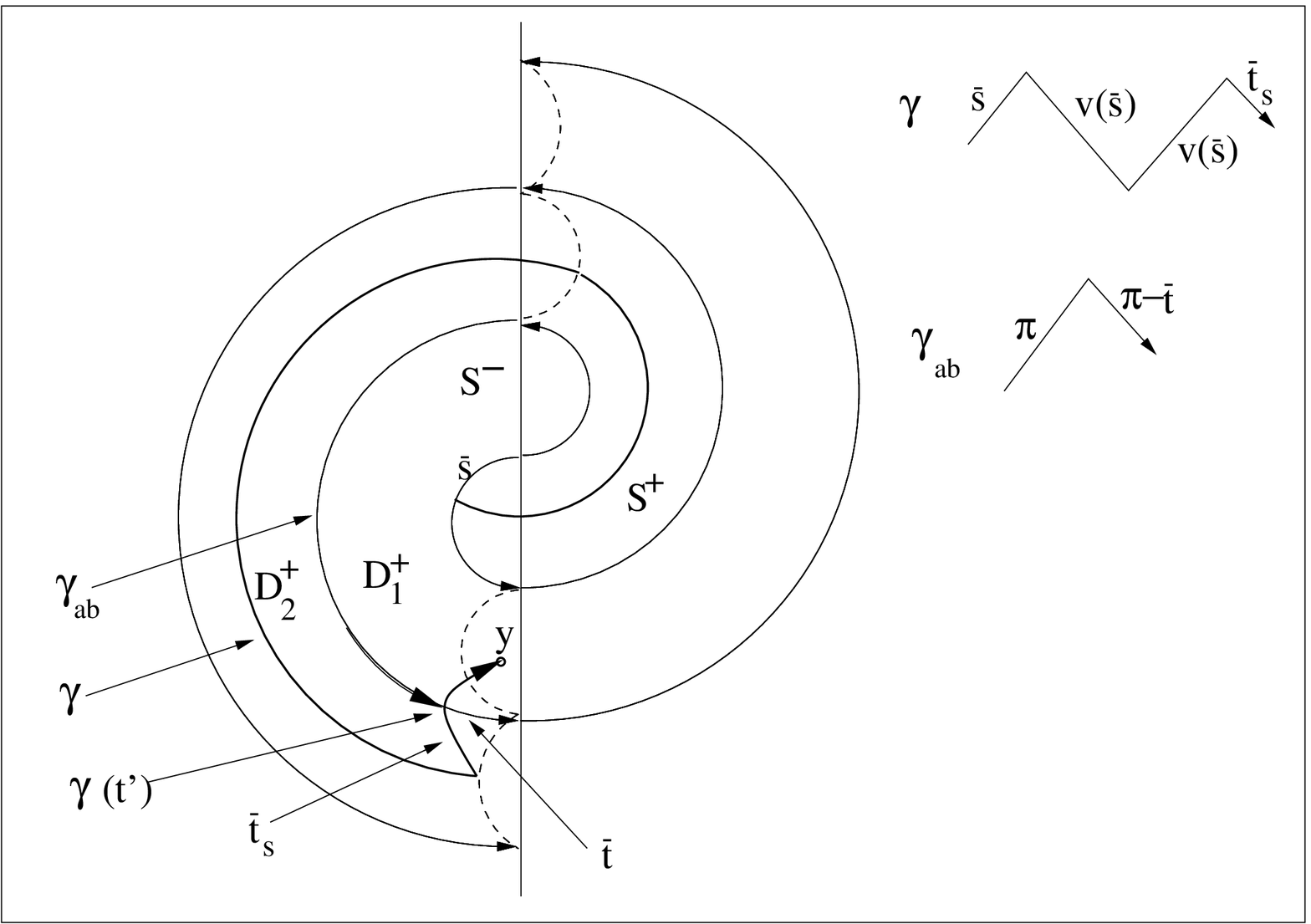}{Proof of Lemma \ref{l-snake-op}. Here to fix the ideas we 
set $\eps=-,~\eps'=+, ~k=1$}{11}
Proof of Claim 1: it is clear that there exists a neighborhood $U$
of $t'$ so that $\ga|_U$ is an integral curve of $X_S^{-\eps}$. Thanks
to Remark~\ref{r-rel-conc} and to the argument above, $\ga|_U$ intersects 
$Int(L^+_k)$ transversally (see figure \ref{f-lemmi})
in such a way that $\ga$, run backward in time, goes from $D^+_k$ to
$D^+_{k+1}$. Claim 1 is proved.\quadv

Now, by definition of $t'$, $\g(t')\in{\cal A}^\eps|_{[0,T_{eq}]}$.
Let $\ga_{ab}$ be the restriction of ${\cal{A}}^{\eps}$ between
$y_0$ and $\ga(t')$. Consider $\tilde{\ga}$, the concatenation of 
$\ga_{ab}$ and $\ga|_{(t',T]}$. The conclusion of Lemma~\ref{l-snake-op}
will follow if one can show that the time duration $T'$ of $\tilde{\ga}$
is less than $T$, the time duration of $\ga$. This, in turn, amounts to
show that $T'$, the time duration of $\ga_{ab}$ is less than $t'$,
the time duration of $\ga|_{[0,t']}$. This is the object of the next 
Claim.\ffoot{non e' che si capisca proprio tanto questa 
dimostrazione}\\\\
{\bf Claim 2}
{\it With the notations above, we have $T'<t'$.}\\\\
Proof of Claim 2: The trajectory $\ga$, run backward in time from $t'$, is 
an 
$X^{-\eps'}_S$-integral curve until it hits a (SC)
in some $D^{\eps}_L\in \NH^{\eps'}$, for some integer $L\geq k+1$, at a 
point $C^{\eps'}_L(\bar{s})$,
$\bar{s}\in ]0,\pi[$. 
One can easily conclude that the only possibility is $L=k+1$.
By Claim 1, a time optimal trajectory can 
pass (backward in time) from $S^{\eps}$ to $S^{-\eps}$ only by going 
down, i.e. by 
passing 
from some $D^{\eps'}_k$ to $D^{\eps'}_{k+1}$. Therefore, by an elementary
counting argument, one gets\ffoot{qui c'era scritto $\geq$, ma per me e' 
$=$. INOLTRE $s'=s$.!!!!!!!!!!}
$$
t'=\bar{s}+t_{\bar{s}}+(k+1)v(\bar{s}),
$$
where $t_{\bar{s}}$
is the time needed to go from $\ga(t')$ to $C^+_{k+1}(\bar{s})$.
On the other hand, 
$$
T'=(K+1)\pi-\tilde{t},
$$
where $\tilde{t}\in ]0,\pi[$ is the time needed to ${\cal A}^\eps$ to go 
from $\ga(t')$ to $P^+_{k+1}$. Since $v(\bar{s})>\pi$, then $T'<t'$.
The proof of Lemma \ref{l-snake-op} is finished. \quadp
\brem\llabel{ab-op}
Coupled with the proof of Claim 2, a simple continuity argument implies
that ${\cal{A}}^+$ and ${\cal{A}}^-$ are time optimal trajectories in the top
hemisphere.
\erem

Gathering all the information on time optimal trajectories,
we are now able to compute $N_S(\alpha)$.

\ffoot{QUESTO E" DA VERIFICARE}
\bp\llabel{p00}
For $\alpha\in ]0,\pi/4[$, we have
\begin{equation}\llabel{cor000}
N_S(\alpha):=2\left[\frac{\pi}{8\alpha}\right]-
\left[2\left[\frac{\pi}{8\alpha}\right]-\frac{\pi}{4\alpha}\right].
\end{equation}
\ep
{\bf Proof of Proposition~\ref{p00}}
Let $y\in S^+$ and $\ga$, a time-optimal 
trajectory connecting $y_0$ to $y$. The point $y$ belongs to some $D^+_k$, 
$k\leq N_0$ and, by Lemma~\ref{l-snake-op}, $\ga$ remains in $S^+$.
Since the function $v$ takes values in $[\pi,\pi+\pi/2]$, it is easy to see, 
from Eqs. (\ref{setEX-1}) and (\ref{setEX-2}) and Remark~\ref{ab-op} that $\ga$,
run backward in time, will go through the ordered sequence of regions
$D^+_k$, $D^-_{k-1}$, $D^+_{k-2}$, etc, until hitting one of the two curves
$\ga_{op}^+$ or $\ga_{op}^-$. Moreover, in each of the region $D^{\eps}_l$, 
$\ga$ will switch {\it exactly} once, thanks to Proposition~\ref{p-switchings}.
Therefore, the number of times, where an optimal trajectory 
$\ga$ starting at $y_0$ switches, is exactly equal to the number 
of times $\ga$ crosses the subset of the meridian contained in $\NH$.
The same conclusion holds
for points belonging to $S^-$, ${\cal{A}}^+$ and ${\cal{A}}^-$.
By a systematic examination of all the possibles cases, we end up with 
(\ref{cor000}). Note that $N_S(\alpha)$ is the number of switchings for a time
optimal trajectory ending on the equator.
\quadp

\subsection{Geometric Remarks}\llabel{pendolo}
\subsubsection{Relations with the Linear Pendulum}
In a fixed \neigh of the north pole, the control system 
on the sphere 
\r{giu} behaves, when $\al>0$ is small enough,
as a controlled linear pendulum.
More precisely, let us consider the stereographic projection of the sphere 
from the south pole $(0,0,-1)$ on $V$, the tangent plane to the sphere at 
the north pole.
If $y_1,y_2,y_3$ are the coordinates of the three dimensional Euclidean 
space where the sphere is embedded, a system of coordinates on $V$ 
is $(y_1,y_2)$, see Fig. \ref{f-pendolo}. Let $x_0^+$ and $x_0^-$ be the 
projections of the equilibrium points of $X^+_S$ and  $X^-_S$ in $\NH$.
\ppotR{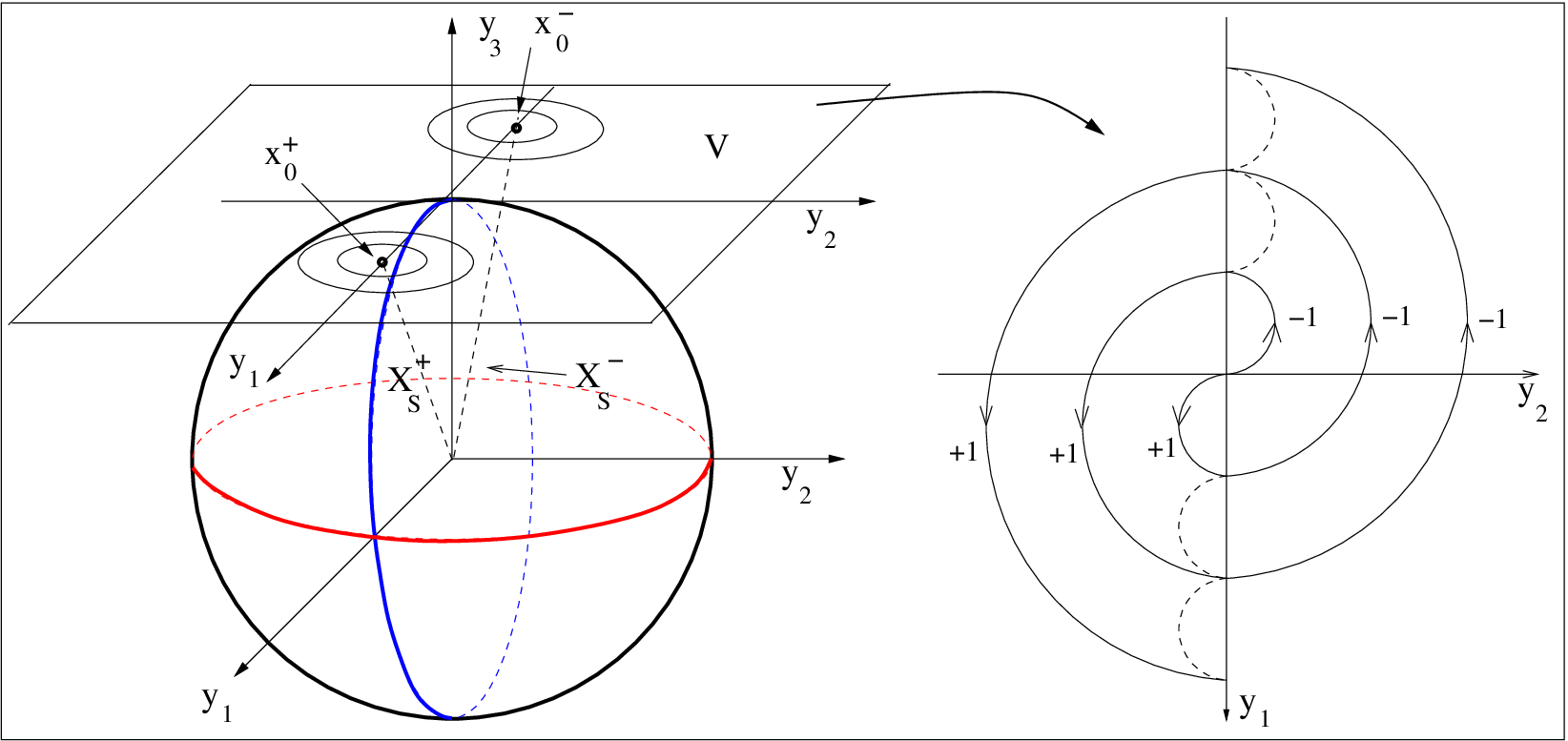}{Stereographic projection and synthesis of the linear 
pendulum}{14.6}
An alternative way of parametrizing this problem (instead of fixing the 
radius of the sphere and varying the axes of rotations) consists of fixing 
the  points $x_0^+=(1,0)^T,$ $x_0^-=(-1,0)^T$ 
and varying the radius $r$ of the sphere. The relation between $\al$ and 
$r$ is $\tan(\al)=1/r.$
The range $\al\in]0,\pi/4[$ becomes $r\in]1,\infty[$ and 
$\al\to0$ corresponds to $r\to\infty$.
In $V$, fix a ball $B(0,r_0)$ of radius $r_0>0$ 
centered in the origin, and consider the stereographic projection of the 
integral curves 
of $X^+_0$ and $X^-_0$. For $r\to\infty$, they become circles centered at 
the  points $x_0^\pm$. Then, one easily sees that, in $B(0,r_0)$, 
the limit system (and the associated  
synthesis) corresponds to a controlled linear pendulum (with the 
associated 
synthesis) of equation:
$\dot y_1=-y_2,$ 
$\dot y_2=y_1-u,$ $|u|\leq1.$
Notice that $\lim_{\al\to\infty}v(s)=\pi$, that is exactly the 
time duration of 
interior bang arcs for the linear pendulum.

\subsubsection{ The time optimal problem on $SU(2)$}

The optimal control problem on $\underline{NH}$ is the projection 
(by a Hopf fibration) of an optimal control problem on $SO(3)$.
Similarly, the corresponding problem on the whole sphere $S^2$ is the projection 
(by an appropriate Hopf fibration) of an optimal control problem on $SU(2)$.
Indeed, $SU(2)$ is the universal (double) covering of $SO(3)$ and they have
the same Lie algebra $so(3)$. 
{The existence of that double covering justifies, by a factor two, 
the difference between our bound and the bound
\r{agr0},  on the maximal number of 
switchings for the control problem on $SO(3)$.} 
Indeed, the index theory developed by Agrachev and Gamkrelidze in 
\cite{agra-sympl-x,agra-sympl} provides a bound on 
the number of switchings by proving that a certain extremal is not optimal
because it loses {\sl local optimality} working at the Lie algebraic level.
This is why the upper bound in \r{agr0} corresponds (essentially) to a control
problem on $SU(2)$, and thus, after projection, on a 
control problem on the whole sphere $S^2$, and not just on $\underline{NH}$.
{The other factor two, of the difference between our bound and the 
bound given in \r{agr0}, comes from the fact that in \cite{agra-sympl-x} 
the index of the second variation was estimated up to an additive factor 1 
(see \cite{agra-sympl-x}, 
p.275).} 

\section{Conclusion and Open Problems}
\llabel{s-5}
In the previous Section, we derived a set of properties of the optimal 
synthesis that were sufficient to compute the maximum number of switchings 
of a time optimal trajectory joining $y_0$ to any point of the north 
hemisphere. This enabled us to provide a precise estimate for $N(\al)$, 
$\al\in]0,\pi/4[$.
However, the following questions remain unsolved:
\bd
\i[Question 1] are all the extremal trajectories \r{setEX-1} 
and \r{setEX-2} optimal in the north hemisphere? 
\ed
The answer to this question depends on the answer to the next question:
\bd
\i[Question 1'] in the north hemisphere, are 
the switching curves $C^\eps_k(s)$, $s\in]0,\pi[$, locally optimal? (The 
points $s=0,s=\pi$ are not included since we already know that the two 
abnormal extremal ${\cal A}^\pm$ are optimal in $\NH$.)
\ed
Roughly speaking we say that a switching curve is locally 
optimal if it never ``reflects'' the
trajectories. 
More precisely, we have the following definition (clarified 
by Fig.\ref{f-op1bis}). 
\bdeff
\llabel{d-opC}
Consider a smooth switching curve $C$ between two smooth vector field
$Y_1$ and $Y_2$ on a smooth two dimensional manifold.
Let $C(s)$ be a
smooth parametrization of $C$.
We say that $C$
is \underline{locally optimal} if, for every $s\in Dom(C)$, we have
\bqn
\dot C(s)\neq\al_1 Y_1(C(s))+\al_2 Y_2(C(s)), \mbox{ for every
}\al_1,\al_2\mbox{ s.t. }\al_1\al_2\geq0.
\eqnl{conj}
\edeff
\ppotR{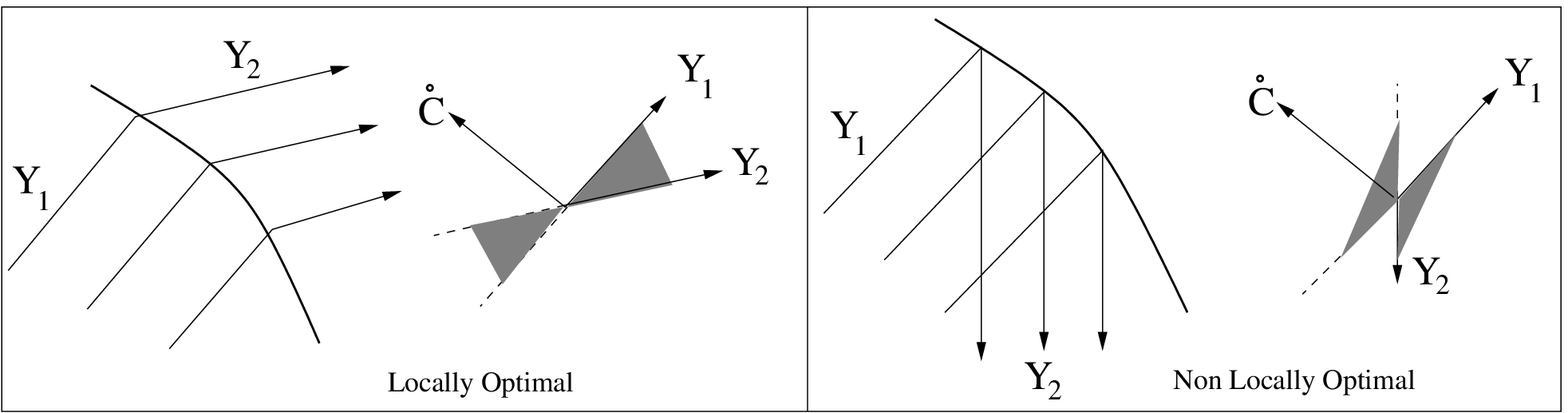}{Definition \ref{d-opC}}{13}
{The points of a switching curve on which relation \r{conj} is not 
satisfied are usually called ``conjugate points''.}
\brem
Notice that, if all the switching curves are locally 
optimal in 
the north  hemisphere, it follows that the set of extremals \r{setEX-1}
and \r{setEX-2} (restricted to $\UNH$) is an 
optimal synthesis for the problem 
\r{giu} on $\R P^2$. In this case, on $\R P^2$, the extremals 
\r{setEX-1}
and \r{setEX-2} lose global  
optimality before losing local optimality. 
\erem
\bd
\i[Question 2] If the answer to Question 1' is yes,
               what about the same question for the optimal control 
               problem on $S^2$?
               \ffoot{(of course in this case the shape of the optimal 
               synthesis 
               will be similar to the one in the north pole except for a 
               \neigh of the south pole).}
               More precisely, one 
               would like to understand  how the extremal trajectories  
               \r{setEX-1} and \r{setEX-2} are 
               going to lose optimality in a \neigh of the south pole 
               (i.e. if the loss of optimality is local or just global).

\i[Question 3] What is the shape of the optimal synthesis in a \neigh of 
the south pole?

\ed
In this section, we present the results of some numerical simulations
which provide some hints regarding the above questions.
More precisely:

\bi 

\i There is strong numerical 
evidence for a positive answer to {\bf Question 1'}. 
This means that the switching curves in the north hemisphere never 
reflect trajectories. In other words, situations like those considered 
in  the proof of Lemma \ref{l-snake-op}, (cf.  Fig.\ref{f-lemmi}) 
are not possible. 

\i As concerns {\bf Question 2}, we conjecture the following:

\bd
\i[C1] 
{The curves $C^\eps_k(s)$, $s\in]0,\pi[$ are locally optimal if and 
only if
$X^+_S(C^\eps_k(0))=\al_1 X^-_S C^\eps_k(0))$ and 
$X^+_S(C^\eps_k(\pi))=\al_2 X^-_S C^\eps_k(\pi))$ with $\al_1,\al_2\geq 0$ 
but not both vanishing.}
\ed
{This condition is verified if and only if:}
$k\leq\left[\frac{\pi-\al}{2\al}\right]-1$, 
which simply follows from Remark \ref{r-rel-conc}.

Set  $\NA:=\left[ \frac{\pi}{2\al}\right]$. Analyzing the evolution of the 
minimum time wave front in a \neigh of the
south-pole, it is reasonable to conjecture that:
\bd
\i[C2] For $T\leq (\NA-1)\pi$, the synthesis built above is 
optimal. Every $x\in S^2$ is reached in time 
$T\leq(\NA+1)\pi$.\ffoot{penso 
che nei casi $B$ e $D2$ si abbia solo $N$, ma controllare}
Every optimal trajectory has at most $\NA$ switchings and
there exists an optimal trajectory having  $\NA-1$ switchings.
\ed
In the top of Fig. \ref{f-critica-unica}, the optimal synthesis is 
plotted.


\i Regarding {\bf Question 3}, numerical simulations suggest that
the shape of the optimal synthesis for time $T>(\NA-1)\pi$, 
depends on the remainder:
$$
r:=\pi-2\al \NA=\pi-2\al \left[ \frac{\pi}{2\al}       
\right].
$$
Notice that $r$ belongs to the interval $[0,2\al[$.
More precisely, we conjecture the following:
\bd
\i[C3] For $\al\in]0,\pi/4[$, there exist two positive numbers $\al_1$ 
and $\al_2$ such that $0<\al_1<\al<\al_2<2\al$ and: 
\bd
\i{}\underline{CASE A: $r\in]\al_2,2\al[$}.
The switching curve starting at $P^+_{\NA}$ glues to an 
overlap  curve that passes  through the origin (see the bottom of Fig. 
\ref{f-critica-unica}, Case A).

\i{}\underline{CASE B:  $r\in[\al_1,\al_2]$}. An overlap 
curve starts exactly at $P^+_{\NA}$ and
passes through the origin.

\i{}\underline{CASE C:  $r\in]0,\al_1[$}. The situation is 
more complicated and it is
depicted in the bottom of Fig. \ref{f-critica-unica}, Case C.

\ed

{For $r=0$, the situation is the same as in CASE A, but for the 
switching curve starting at $P^+_{\NA-1}$.}

\ed
\ei


\ffoot{\brem
We believe that for the question 1' it is  possible to 
give a mathematical proof. One possibility is to use the index theory 
developed by Agrachev and his coauthors in 
\cite{agra-sympl-x,agra-sympl,agra-sz}.
Notice that thanks to the recent paper \cite{agra-sz}, index Theory 
provides also a sufficient condition for local optimality (for 
trajectories $\con^0$-close to the reference trajectory).
\ffoot{forse bisogna aggiungere 
ancora un articolo della stefani con zezza} On the other side, giving a 
mathematical proof to question 2 (although index theory can again 
help) is harder in our opinion and maybe could be done for small values 
of the parameter $\al$.  
\erem}

\ppotR{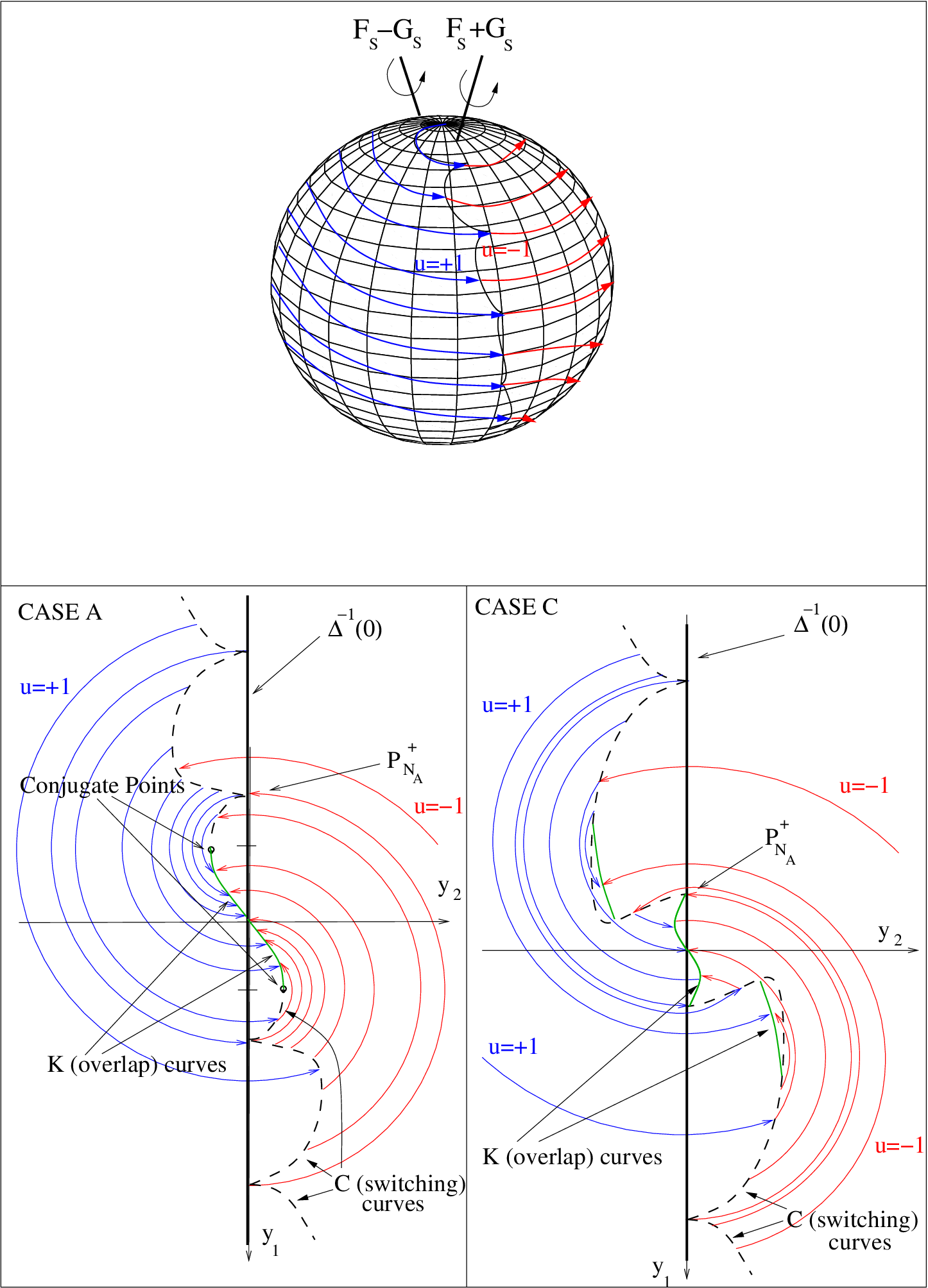}{The Time Optimal Synthesis on the Sphere 
(top) and Optimal Synthesis in a \neigh of the south pole, Case A and Case 
C (bottom)}{13}

\noi
{\bf AKNOWLEDGMENTS} \\
{The authors are grateful to Andrei Agrachev, for 
suggesting the problem and for many geometric hints.  The authors 
would like also to thank Benedetto Piccoli, Gregorio Falqui, Mario 
Sigalotti and Paolo Mason for 
helpful discussions.}



\fine